\documentclass[g5paper,twoside,final,10pt]{article}

\usepackage{amsfonts}
\usepackage{amsmath}
\usepackage{amsthm}
\usepackage{ifthen}
\newboolean{draftversion}

\newcommand{\note}[2][]{\message{(#1)}\ifthenelse{\boolean{draftversion}}%
	{\noindent{\em[#2]}}{}}

\def \P{\mathbb{P}}

\def \E{\mathbb{E}}

\def \eop {\hbox{}\nobreak\hfill
\vrule width 2mm height 2mm depth 0mm
\par \goodbreak \smallskip}

\def \eop {\hbox{}\nobreak\hfill \vrule width 2.0mm height 1.8mm depth 0mm
\par \goodbreak \smallskip}

\numberwithin{equation}{section}

\theoremstyle{plain}
\newtheorem{definition}{Definition}[section]
\newtheorem{theo}{Theorem}[section]
\newtheorem{pro}{Proposition}[section]
\newtheorem{lem}{Lemma}[section]

\newtheorem{ex}{Example}[section]

\def \P{\mathbb{P}}

\def \E{\mathbb{E}}

\def \eop {\hbox{}\nobreak\hfill
\vrule width 2mm height 2mm depth 0mm
\par \goodbreak \smallskip}

\def \P{\mathbb{P}}

\def \E{\mathbb{E}}

\def \eop {\hbox{}\nobreak\hfill
\vrule width 2mm height 2mm depth 0mm
\par \goodbreak \smallskip}

\def \eop {\hbox{}\nobreak\hfill \vrule width 2.0mm height 1.8mm depth 0mm
\par \goodbreak \smallskip}
\numberwithin{equation}{section}

\def \P{\mathbb{P}}

\def \E{\mathbb{E}}

\def \eop {\hbox{}\nobreak\hfill
\vrule width 2mm height 2mm depth 0mm
\par \goodbreak \smallskip}

\def \eop {\hbox{}\nobreak\hfill
\vrule width 1.4mm height 1.4mm depth 0mm
\par \goodbreak \smallskip}

\def \P{\mathbb{P}}

\def \E{\mathbb{E}}

\def \eop {\hbox{}\nobreak\hfill\vrule width 2mm height 2mm depth 0mm
\par \goodbreak \smallskip}

\def \eop {\hbox{}\nobreak\hfill \vrule width 2.0mm height 1.8mm depth 0mm
\par \goodbreak \smallskip}

\title{A maximum principle for relaxed stochastic control of linear SDEs with
application to bond portfolio optimization}
\begin{document}
\author{Daniel Andersson\thanks{Department of
Mathematics, Royal Institute of Technology, S-100 44 Stockholm,
Sweden. danieand@math.kth.se} ~and Boualem
Djehiche\thanks{Department of
Mathematics, Royal Institute of Technology, S-100 44 Stockholm,
Sweden. boualem@math.kth.se}}
\maketitle

\begin{abstract} \noindent We study relaxed stochastic control problems where the state equation is a one dimensional linear stochastic differential equation with random
and unbounded coefficients. The two main results are existence of an optimal relaxed control and necessary conditions for optimality in the form of a relaxed maximum principle.
The main motivation is an optimal bond portfolio problem in a market where there exists a continuum of bonds and the
 portfolio weights are modeled as measure-valued processes on the set of times to maturity.
\bigskip
\bigskip

\textbf{Keywords.} Stochastic control, relaxed control, maximum principle, $\cal H$-function, bond portfolio.

\bigskip
\textbf{AMS subject classification.}  93E20, 60H30, 60H10, 91B28.
\end{abstract}


\section{Introduction}
The objective of this paper is to derive necessary conditions for optimality in relaxed stochastic control problems, i.e.~the control is a measure-valued process,
where the state process is a solution to a one dimensional linear stochastic differential equation (SDE) whose coefficients are random and not necessarily bounded.
This study is motivated by the following optimal bond portfolio problem. Consider a market of non-defaultable bonds, i.e. financial contracts that are bought today and pay a fixed amount at some
 future time, called the maturity time. At each time $t$, the investor is allowed to buy bonds with any time to maturity in $U$, where $U$ is a subset of $\mathbb{R}_+$. Modeling the prices of the bonds as SDEs, we may write down the wealth of the investor as an SDE of the form (see Section \ref{hjmm} below)
\begin{align*}
x_t=x_0+\int_0^tx_s\int_U\left(r_s^0-v_s(u)\Theta_s\right)\mu_s(du)ds+\int_0^tx_s\int_Uv_s(u)\mu_s(du)dB_s,
\end{align*}
where $x_0$ is the investors initial capital and $\mu_t$ is a probability measure on $U$ reflecting the proportion invested in bonds with different maturities. Further, $r_t^0$ is
the short rate, $v_t$ is the integrated volatility process of the bond prices and $\Theta_t$ is the so called market price of risk. The objective of the investor is to choose $\mu_t$ in some optimal way. Interpreting $\mu_t$
as the control process leads to a relaxed control problem where the state process is a linear SDE with random coefficients, and where $r_t^0$ and $\Theta_t$ cannot in general
be assumed to be bounded.
\\

\noindent This motivates us to study
relaxed stochastic control problems where the state equation is a one dimensional linear SDE
 \begin{align*}
 x_t=x_0+\int_0^t\int_Ub(s,x_s,u)\mu_s(du)ds+\int_0^t\int_U\sigma(s,x_s,u)\mu_s(du)dB_s,
 \end{align*}
on some filtered probability space $(\Omega,\mathcal{F},\mathcal{F}_t,\P)$
equipped with a $d$-dimensional Brownian motion $B_t$. The control variable is a process $\mu_t$ taking values in the space of probability measures on the action space $U$.
 $x_0$ denotes the
 initial state, $b$ and $\sigma$ are random coefficients of the form
\begin{align*}
b(t,x,u,\omega)=\upsilon_t(u,\omega)+\phi_t(u,\omega)x,\\
\sigma(t,x,u,\omega)=\chi_t(u,\omega)+\psi_t(u,\omega)x,
\end{align*}
for given stochastic processes $\upsilon$, $\phi$, $\chi$ and $\psi$ taking values in some space of functions on $U$.
The cost functional, which is to be minimized, is of the form
\begin{align*}
 J(\mu_t)=\E\left(\int_0^T\int_Uh(t,x_t,u)\mu_t(du)dt+g(x_T)\right).
 \end{align*}

\medskip
Under the usual assumptions on $b$ and $\sigma$, i.e.~deterministic functions of $(t,x,u)$, Lipschitz continuous and with linear growth in $x$, a maximum principle for stochastic (strict) control problems where the
control enters the diffusion coefficient was established in Peng \cite{peng}. An extension to relaxed control problems is given in Bahlali et al.~\cite{bd}. We refer to Cadenillas and Karatzas \cite{candenillas} for a result on stochastic (strict) control, in terms of a maximum principle for linear SDEs with bounded random coefficients, under integrability conditions on the control.\\

 This paper contains two main results. The first one, Theorem \ref{existence}, establishes existence of an optimal relaxed control which is derived using a similar scheme as in Ma and Yong \cite{mayong}. The main tools in the proof are tightness and Skorohod's Selection Theorem. The second main result, Theorem \ref{pontryagin}, suggests necessary conditions for optimality that are given in form of a relaxed maximum principle. The proof is based on the Chattering Lemma, which gives a sequence of ordinary (strict) controls that approximates the relaxed control. The proof of the maximum principle is based on Zhou's maximum principle \cite{zhou} for near optimal strict controls, and stability properties of the state- and adjoint processes with respect to the control. Note that the relaxed control problems studied e.g. in Bahlali et al.~\cite{bd} is different to ours, in that they relax the corresponding infinitesimal generator of the state process, which leads to a martingale problem for which the state process is driven by an orthogonal martingale measure. In our setting the driving martingale measure $\mu_t(du)dB_t$ is however not orthogonal.\\

 The paper is organized as follows. In Section 2, we formulate the relaxed control problem for our linear SDEs. In Section \ref{existenceresult} we prove existence of an optimal control, while in Section \ref{necessary}, necessary conditions for optimality are given in form of a relaxed maximum principle. In Section 5, we apply these results to formulate a maximum principle for an optimal bond portfolio problem. Finally, to make the exposition simple, all the proofs and technical details are collected in Section 6.
 
\section{Formulation of the problem}
Consider a one dimensional SDE on some probability space $(\Omega,\mathcal{F},\mathcal{F}_t,\P)$:
\begin{align}
x_t=x_0+\int_0^tb(s,x_s,u_s)ds+\int_0^t\sigma(s,x_s,u_s)dB_s,  \label{sde11}
\end{align}
where $x_0\in\mathbb{R}$ is the initial state, $B$ is a $d$-dimensional Brownian motion and $\mathcal{F}_t$ is the $\P$-augmentation of the natural filtration $\mathcal{F}_t^B$ defined by
\begin{align*}
\mathcal{F}_t^B=\sigma\left(B_s,s\in[0,t]\right)~\textrm{for all}~t\in[0,T].
\end{align*} 
 Furthermore, the coefficients are given by
\begin{subequations} \label{linearcoeff}
\begin{align}
b(t,x,u,\omega)=\upsilon_t(u,\omega)+\phi_t(u,\omega)x \label{first} \\
\sigma(t,x,u,\omega)=\chi_t(u,\omega)+\psi_t(u,\omega)x, \label{second}
\end{align}
\end{subequations}
where $\upsilon:\mathbb{R}^+\times U\times\Omega\mapsto \mathbb{R}$, $\phi: \mathbb{R}^+\times U\times\Omega
\mapsto \mathbb{R}$, $\chi: \mathbb{R}^+\times U\times\Omega\mapsto \mathbb{R}^d$ and
$\psi:\mathbb{R}^+\times U\times\Omega\mapsto \mathbb{R}^d$ are $\mathcal{F}_t$-adapted
processes. For each $t\in[0,T]$ the control $u_t$ is in the action space $U$, a compact set in $\mathbb{R}^n$. Let $\mathcal{U}$ denote the class of admissible controls, i.e.~$\mathcal{F}_t$-adapted processes with values in $U$. The cost functional is given by
\begin{align}
J(u)=\E\left(\int_0^Th(t,x_t,u_t)dt+g(x_T)\right), \label{costfunctional1}
\end{align}
and the objective is to minimize $J$ over the set of admissible controls. A control $u^*$ is called optimal if it satisfies $J(u^*)=\inf\{J(u);u\in\mathcal{U}\}$. If also $u^*\in
\mathcal{U}$, it is called a strict optimal control.

\medskip
We make the following assumptions regarding the state equation (\ref{sde11}), (\ref{linearcoeff}) and cost functional (\ref{costfunctional1}).
\begin{itemize}
\item[(A.1)] $\varphi_t(u,\omega)$ is continuous in $(t,u)$, where $\varphi$ stands for one of the processes $\upsilon, \phi, \chi, \psi$.
\item[(A.2)] $\psi$ is uniformly bounded in $\mathbb{R}^+\times U\times\Omega$.
\item[(A.3)] For any $k\in(-\infty,\infty)$ it holds that 
\begin{align*}
\E\left(\exp\left[k\int_0^T\phi_t(u)dt\right]\right)<\infty, ~\mbox{for all}~ u\in U.
\end{align*}
\item[(A.4)] For any $k>0$, it holds that 
\begin{align*} 
\E\left(\int_0^T|\upsilon_t(u)|^kdt\right)<\infty\,\,\, \mbox{and}\,\,\, \E\left(\int_0^T|\chi_t(u)|^kdt\right)<\infty,
\end{align*}
for all $u\in U$.
\item[(A.5)] The functions $g$ and $h$ are twice continuously differentiable in $x$.~The function $g$ and its first and second derivative are bounded and Lipschitz continuous in $x$. The function $h$ and its first and
second derivative are bounded, continuous in $u$ and Lipschitz continuous in $x$.
 \end{itemize}

\medskip
\noindent Throughout the rest  of the paper we will not specify that properties hold $\P$-a.s. when it is clear from the context. We denote for any process $\varphi_t$, \begin{align*}
\vert\varphi\vert^{*,p}_T=
\sup_{t\in[0,T]}|\varphi_t|^p.
\end{align*}

\medskip\noindent This kind of control problems is often formulated in the so-called relaxed form, due to the fact that a strict optimal control may fail to exist (see e.g. Bahlali et al.~\cite{bd} for a discussion). Instead one embeds the strict
controls in a wider class of controls that takes values in probability measures on $U$ rather than on $U$ itself. Also, a solution to a relaxed control problem is a weak one, i.e.~the probability space, equipped with the a priori given stochastic processes, is part of the solution.\\

\medskip\noindent Let $\mathcal{P}( U)$ be the space of probability measures on $U$. If
$\mu_t(du)$ is a stochastic process taking values in $\mathcal{P}(U)$, we denote by $\mathcal{L}([0,T],U)$ the space of the (Radon) measure-valued processes
$d\lambda_t(u)=\mu_t(du)dt$. If a probability space $(\Omega, \mathcal{F}, \P)$ is given, then we denote $\boldsymbol{M}(\Omega)$ the space of all $\mathcal{F}_t$-adapted processes $\mu_t(du)$ taking values in $\mathcal{P}( U)$. Further, we denote by $\boldsymbol{L}(\Omega)$ the space of all $\mathcal{L}([0,T],U)$-valued $\mathcal{F}_t$-adapted processes. It can be shown that there is a one-to-one correspondence between  $\boldsymbol{M}(\Omega)$ and $\boldsymbol{L}(\Omega)$, and that $\mathcal{L}([0,T],U)$ is a compact metric space. For further discussion, see Ma and Yong \cite{mayong}. \\

\medskip\noindent Throughout we denote $f(\mu_t)=\int_Uf(u)\mu_t(du)$, for any continuous function $f$.
By expanding the set of controls from $\mathcal{U}$ to $\boldsymbol{M}$,
  the state equation is defined as
 \begin{align}
 x_t=x_0+\int_0^tb(s,x_s,\mu_s)ds+\int_0^t\sigma(s,x_s,\mu_s)dB_s. \label{sde1}
 \end{align}

\begin{definition} \label{defrelaxedcontrol}
A relaxed control is the term $\mathcal{A}=(\Omega, \mathcal{F}, \mathcal{F}_t, \P, \boldsymbol{B}_t, \mu_t, x_t) $, where
\begin{itemize}
\item[$(i)$] $(\Omega, \mathcal{F}, \mathcal{F}_t, \P)$  is a filtered probability
space;
\item[$(ii)$] $\boldsymbol{B}_t= (B_t,\upsilon_t, \phi_t, \chi_t, \psi_t)$, in which $B_t$ is an $\mathcal{F}_t$-Brownian motion and
$\upsilon_t, \phi_t, \chi_t, \psi_t$ are $\mathcal{F}_t$-adapted stochastic processes satisfying \textnormal{(A.1)-(A.4)};
\item[$(iii)$] $\mu_t\in
\boldsymbol{M}(\Omega)$;
\item[$(iv)$] $x_t$ is $\mathcal{F}_t$-adapted and satisfies \textnormal{(\ref{sde1})}.
\end{itemize}
\end{definition}

\medskip
\noindent We denote by $\mathcal{U}^R$ the set of all relaxed controls. The cost functional corresponding to the control $\mathcal{A}$
 is defined as
\begin{align}
J(\mathcal{A})=\E\left(\int_0^Th(t,x_t,\mu_t)dt+g(x_T)\right), \label{costfunctional}
\end{align}
and a relaxed control $\mathcal{A}^*$ is optimal if $J(\mathcal{A}^*)=\inf\{J(\mathcal{A});\mathcal{A}\in\mathcal{U}^R\}$. It is well known that $\mathcal{U}$ may be embedded into
$\mathcal{U}^R$, since any strict ($U$-valued) control process $u_t$ can be represented as a relaxed control by setting $\mu_t(du)=\delta_{u_t}(du)$. Moreover the so-called
Chattering Lemma, stated in Section \ref{necessary}, tells us that any relaxed control is a weak limit of a sequence of strict controls.


\section{Existence of an optimal relaxed control}\label{existenceresult}
In this section we shall establish the existence of an optimal relaxed control. To achieve this, we construct a minimizing sequence of controls $\mathcal{A}^{(k)}\in\mathcal{U}^R$ for the cost functional $J$, i.e.
\begin{align*}
\inf\{J(\mathcal{A}),\mathcal{A}\in\mathcal{U}^R\}=\lim_{k\to\infty}J(\mathcal{A}^{(k)}),
\end{align*}
and show that a limit $\mathcal{A}$ exists and fulfills $(i)-(iv)$ in Definition \ref{defrelaxedcontrol}.
This will be carried out i several steps described in Lemmas \ref{contmeasure}-\ref{weakconv} below, (cf. the scheme suggested in Ma and Yong \cite{mayong}). The main tools are tightness of the processes and Skorohod's
 Selection Theorem.  To make the exposition simple, we will consider a simpler form the cost functional $J$, by letting $h=0$. The proofs can be modified so that the results hold without this restriction. This is
done by adding $\int_0^T\int_Uh(x_t,u)\mu_t(du)dt$ as an "extra state" (cf. Yong and Zhou \cite{yongzhou}).

\medskip\noindent
Given a relaxed control $\mathcal{A}=(\Omega,\mathcal{F},\mathcal{F}_t,\P,\boldsymbol{B}_t,\mu_t,x_t)$, there exists a unique strong solution to the equation given by
(\ref{linearcoeff}) and (\ref{sde1}). Moreover, its explicit form is
\begin{align} \label{explicit}
x_t=z_t\left(x_0+\int_0^t\frac{\upsilon_s(\mu_s)-\psi_s(\mu_s)\chi_s(\mu_s)}{z_s}ds+\int_0^t\frac{\chi_s(\mu_s)}{z_s}dB_s\right),
\end{align}
where
\begin{align}\label{explicit2}
z_t=\exp\left(\int_0^t\left(\phi_s(\mu_s)-\frac{1}{2}\psi_s^2(\mu_s)\right)ds+\int_0^t\psi_s(\mu_s)dB_s\right).
\end{align}
This can be verified by applying Ito's formula on (\ref{explicit})-(\ref{explicit2}). Moreover, with the assumptions (A.1)-(A.4) we can prove by standard methods that
$x_t$ has the following properties:
For any $p\geq1$ we have
\begin{align}
\E|x|_T^{*,p} < \infty, \label{momentsde}
\end{align}
and there exists a constant $K>0$ such that
\begin{align}
\E(x_t-x_s)^4\leq K|t-s|^2, \label{kolmogorovcondition}
\end{align}
for all $s,t \in [0,T]$.

\medskip
\begin{lem}\label{contmeasure} Given a relaxed control
$\mathcal{A}=(\Omega,\mathcal{F},\P,\mathcal{F}_t,\boldsymbol{B}_t,\mu_t, x_t)$, there exists a sequence $\mu_t^{(k)} \in \boldsymbol{M}(\Omega)$ such that for each $k$,
 the path $\mu_t^{(k)}(A)$ is
continuous for all Borel sets $A$; and
\begin{align}
 \E|x^{(k)}-x|_T^{*,2}\to 0, ~ \emph{as} ~ k\to \infty, \label{approx}
 \end{align}
where $x_t^{(k)}$ is the solution to \textnormal{(\ref{sde1})} with respect to $\mu_t^{(k)}$. Consequently, setting \\
$\mathcal{A}^{(k)}=
(\Omega,\mathcal{F},\P,\mathcal{F}_t,\boldsymbol{B}_t,\mu_t^{(k)}, x_t^{(k)})$, we get
\begin{align*}
J(\mathcal{A}^{(k)})\to J(\mathcal{A}), ~ \emph{as} ~ k\to\infty.
\end{align*}
\end{lem}

\medskip
\noindent Let $\mathcal{A}^{(k)}=(\Omega^{(k)},\mathcal{F}^{(k)},\P^{(k)},\mathcal{F}_t^{(k)},\boldsymbol{B}_t^{(k)},\mu_t^{(k)},x_t
^{(k)})$ be the minimizing sequence, i.e.
\begin{align*}
J(\mathcal{A}^{(k)})\to
\inf\{J(\mathcal{A});\mathcal{A}\in\mathcal{U}^R\}, ~ \textrm{as} ~ k\to \infty. \end{align*}
For this sequence we have the following important property.

\begin{lem}\label{tight}
Denote ~$d\lambda_t^{(k)}(u)=\mu_t^{(k)}(du)dt$. Then the sequence
$(\boldsymbol{B}_t^{(k)},\lambda_t^{(k)},x_t^{(k)})$ is tight in
$\left(C([0,T])\times C([0,T]\times U)^4\right)\times\mathcal{L}([0,T]\times U)\times C([0,T])$.
 \end{lem}

\medskip
 \noindent By Lemma \ref{contmeasure} we may assume that $\mu_t^{(k)}(A)$ has continuous paths for each Borel set $A$. By tightness and the Skorohod's Selection Theorem there exists a probability space
$(\hat{\Omega},\hat{\mathcal{F}},\hat{\P})$, on which is defined a sequence of processes $(\hat{\boldsymbol{B}}_t^{(k)},\hat{\lambda}_t^{(k)},\hat{x}_t^{(k)})$
identical in law to $(\boldsymbol{B}_t^{(k)},\lambda_t^{(k)},x_t^{(k)})$ and converging $\hat{\P}$-a.s. to $(\hat{\boldsymbol{B}}_t,\hat{\lambda}_t,\hat{x}_t)$. Moreover, by Lemma
2.1 in Ma and Yong \cite{mayong} the processes $\hat{\mu}_t^{(k)}$ corresponding to $\hat{\lambda}_t^{(k)}$ have the same law as $\mu_t^{(k)}$ since $\mu_t^{(k)}$ has continuous paths.
We drop the ``$\hat{~~}$'' in the following.

\begin{lem}\label{weakconv}
Let $b$ and $\sigma$ be the processes defined by \textnormal{(\ref{linearcoeff})}, then
\begin{align*} b(t,x_t,\mu_t^{(k)})\stackrel{w}\longrightarrow b(t,x_t,\mu_t),\\
  \sigma(t,x_t,\mu_t^{(k)})\stackrel{w}\longrightarrow \sigma(t,x_t,\mu_t),
 \end{align*}  in $L^2([0,T]\times\Omega)$, as $ k\to
\infty$.
\end{lem}

\noindent\textbf{Proof.}
See Ma and Yong \cite{mayong}, Lemma 3.3.
\eop

\medskip
\noindent Again, using Skorohod's Selection Theorem, there exists a limit $\mathcal{A}$ of the minimizing
 sequence $\mathcal{A}^{(k)}$ which satisfies $(i)-(iii)$ in Definition \ref{defrelaxedcontrol}. By Lemma \ref{weakconv} and the fact that $(\boldsymbol{B}_t^{(k)},\lambda_t^{(k)},x_t^{(k)})\to (\boldsymbol{B}_t,\lambda_t,x_t)$~$\P$-a.s., as $k\to\infty$, we prove in the next theorem that $(iv)$ also holds:

\begin{theo}\label{existence}
The limit ~$x_t$ of $x_t^{(k)}$ satisfies
\begin{align*}
x_t=x_0+\int_0^tb(s,x_s,\mu_s)ds+\int_0^t\sigma(s,x_s,\mu_s)dB_s.
\end{align*}
Therefore, $\mathcal{A}=(\Omega,\mathcal{F},\P,\mathcal{F}_t,\boldsymbol{B}_t,\mu_t,x_t)$ is an optimal relaxed control.
\end{theo}


\section{A relaxed maximum principle} \label{necessary}
We start by proving the so-called Chattering Lemma, which states that any relaxed control may be approximated by strict controls.
 \begin{lem}\label{chattering}$(The~ Chattering~ Lemma)$ ~\\
 Let $\mathcal{A}=
(\Omega, \mathcal{F}, \mathcal{F}_t, \P, \boldsymbol{B}_t, \mu_t, x_t) $ be a relaxed control.
  Then there exist relaxed controls
$\mathcal{A}^{(k)}=
  (\Omega, \mathcal{F}, \mathcal{F}_t, \P, \boldsymbol{B}_t, \delta_{u_t^{(k)}}, x_t^{(k)})$,
 $\hat{\mathcal{A}}^{(k)}=(\Hat{\Omega},\Hat{\mathcal{F}},\hat{\mathcal{F}_t},\hat{\P},\hat{\boldsymbol{B}}_t,\delta_{\hat{u}^{(k)}_t}, \hat{x}_t^{(k)})$ and
 $\hat{\mathcal{A}}=(\Hat{\Omega},\Hat{\mathcal{F}},\hat{\mathcal{F}_t},\hat{\P},\hat{\boldsymbol{B}}_t,\hat{\mu}_t, \hat{x}_t)$, where  $u_t^{(k)}$ and $\hat{u}_t^{(k)}$ are sequences of $U$-valued progressively measurable processes defined on
 $(\Omega,\mathcal{F},\P)$ and $(\hat{\Omega},\hat{\mathcal{F}},\hat{\P})$ respectively, such that
 $\delta_{u_t^{(k)}}$ and $\delta_{\hat{u}_t^{(k)}}$ as well as $\mu_t$ and $\hat{\mu}_t$ are identical in law and such that
if we denote $d\hat{\lambda}_t^{(k)}(u)=\delta_{\hat{u}_t^{(k)}}(du)dt$ and $d\hat{\lambda}_t(u)=\hat{\mu}_t(du)dt$, then
\begin{align}
\hat{\lambda}_t^{(k)}\to\hat{\lambda}_t,
\end{align}
as $k\to\infty$ $\hat{\P}$-a.s. in $\mathcal{L}([0,T]\times U)$. Moreover,
 $x_t^{(k)}$ and $\hat{x}_t^{(k)}$ as well as $\hat{x}_t$ and $x_t$
 are identical in law and
 \begin{align}
 \hat{\E}|\hat{x}^{(k)}-\hat{x}|_T^{*,2}\to 0, \label{chattering1}
 \end{align}
 as $k\to\infty$.
\end{lem}
\medskip
\noindent With the definitions in the Chattering Lemma we thus have
 \begin{align}\label{convergence of J}
 J(\mathcal{A}^{(k)})\to J(\mathcal{A}),
 \end{align}
 as $k\to\infty$, and consequently
 \begin{align}
 \inf\{J(\mathcal{A});\mathcal{A}\in\mathcal{U}^R\}=\inf\{J(u);u\in\mathcal{U}\}. \label{strictrelaxed}
 \end{align}
 The latter equality motivates the use of relaxed control even when one is only concerned with strict controls. The strict and the relaxed problems have the same optimal value.
  However, this optimum may not be reached by a strict control but
 only with a measure-valued one which in turn can be approximated by strict controls.\\

 \noindent By the Chattering Lemma we can assume that any relaxed control has an approximating sequence defined on the same probability space. Let $\hat{\mathcal{A}}=\\(\Omega,\mathcal{F},\P,\mathcal{F}_t,\boldsymbol{B}_t,\hat{\mu}_t, \hat{x}_t)$
be an optimal relaxed control. From now on we let this filtered probability space be fixed and vary only the control measures and corresponding state processes.

\subsection{Adjoint processes}
We recall the first- and second order adjoint processes for the state process (\ref{sde11})-(\ref{linearcoeff}). These are two pairs of processes $(p,q)$ and $(P,Q)$ with values
 in $\mathbb{R}\times \mathbb{R}^{d}$ defined for any strict control $u\in\mathcal{U}$.
We denote by $f_x$ and $f_{xx}$ the first and second derivative, respectively, with respect to $x$ of the
function $f$, where $f$ stands for either $g$ or $h$. Then $(p,q)$ and $(P,Q)$ are given by
\begin{align}
\left\{ \begin{array}{ll}
dp_t=&-\Big(\phi_t(u_t)p_t+\psi_t(u_t)q_t+h_x(x_t,u_t)\Big)dt+q_tdB_t\\
p_T=&g_x(x_T) \label{adjoint11}
\end{array} \right.
\end{align}
\begin{align}
\left\{ \begin{array}{ll}
dP_t=&-\Big(2\phi_t(u_t)P_t+\psi_t(u_t)P_t\psi_t(u_t)+2\psi_t(u_t)Q_t \\
&+h_{xx}(x_t,u_t)\Big)dt+Q_tdB_t\\
P_T=&g_{xx}(x_T) \label{adjoint12}
\end{array} \right.
\end{align}
The second order adjoint process (\ref{adjoint12}) appears when the control affects the uncertainty (noise part) of the system, i.e. when the diffusion coefficient depends explicitly on $u$ (cf. Peng \cite{peng}).
Also, note that the reason for the extra components $q$ and $Q$ is to make it possible to find adapted solutions to the backward SDEs (see Ma and Yong \cite{mayong2} for further discussion). Next, we introduce the Hamiltonian of the system:
\begin{align*}
H(t,x,a,p,q)=-h(x,a)-p\Big(\upsilon_t(a)+\phi_t(a)x\Big)-q\Big(\chi_t(a)+\psi_t(a)x\Big)
\end{align*}
for $(t,x,a,p,q)\in[0,T]\times\mathbb{R}\times U\times\mathbb{R}\times\mathbb{R}^d$. Further, we define the $\mathcal{H}$-function corresponding to a given strict control $u_t$ and its corresponding state process $x_t$ by
\begin{align*}
 \mathcal{H}^{(x_t,u_t)}(t,x,a)=&H\Big(t,x,a,p_t,q_t-P_t\big(\chi_t(u_t)+\psi_t(u_t)x_t\big)\Big)\\
 &-\frac{1}{2}\Big(\chi_t(a)+\psi_t(a)x\Big)P_t\Big(\chi_t(a)+\psi_t(a)x\Big)
\end{align*}
for $(t,x,a)\in [0,T]\times \mathbb{R}\times U$, where $p_t,P_t$ and $q_t$ are determined by the adjoint
equations (\ref{adjoint11}) and (\ref{adjoint12}). The next proposition relates the $\mathcal{H}$-function to the cost functional. It was first proved in Peng \cite{peng} for the case of general SDEs with bounded coefficients.

\begin{pro}\label{varineq}
Let $u\in\mathcal{U}$ with state process $x_t$ be given. Denote by $u^\theta$ the perturbed control:
\begin{align*}
u^\theta_t=\left\{ \begin{array}{ll}
v& \textrm{for}~ t\in [\tau,\tau+\theta]\\
u_t & \textrm{otherwise}.\end{array}\right.
\end{align*}
Then there exists two pairs of processes $(p_t,q_t),(P_t,Q_t)$ which solve \textnormal{(\ref{adjoint11})-(\ref{adjoint12})}, such that
\begin{align*}
\E|p|_T^{*,p}+\E\int_0^T|q|^pdt<\infty,\\
\E|P|_T^{*,p}+\E\int_0^T|Q|^pdt<\infty,
\end{align*}
for any $p\geq 1$, and such that the following holds.
\begin{align*}
J(\delta_{u^\theta})-J(\delta_u)=\E\int_0^T\left(\mathcal{H}^{(x_t,u_t)}(t,x_t,u_t)-\mathcal{H}^{(x_t,u_t)}(t,x_t,u^\theta_t)\right)dt+o(\theta).
\end{align*}
\end{pro}

\subsection{Necessary conditions for near optimality}

By Proposition \ref{varineq} we can derive the following necessary condition for a `near-optimal' strict control in
terms of the $\mathcal{H}$-function.

\begin{pro}\label{nearoptimality}
Let $u_t$ be a strict control such that
\begin{align*}
J(\delta_{u})\leq J(\hat{\mu})+\epsilon,
\end{align*}
then there exists constants $K>0$ and $\gamma>0$ such that the following inequality holds
\begin{align}
\E\int_0^T\mathcal{H}^{(x_t,u_t)}(t,x_t,u_t)dt
\geq &\underset{a\in \mathcal{U}}{\sup}\E\int_0^T\mathcal{H}^{(x_t,u_t)}(t,x_t,a)dt-
K\epsilon^{\gamma}. \label{Hinequality}
\end{align}
\end{pro}

\subsection{The relaxed maximum principle}

We define the adjoint equations corresponding to a relaxed control $\mu_t$ as in (\ref{adjoint11}) and (\ref{adjoint12}) with $\varphi_t(u_t)$ replaced by $\varphi_t(\mu_t)=\int_U \varphi_t(u)d\mu_t(u)$ where $\varphi_t$ is any of $\upsilon_t,\phi_t,\chi_t,\psi_t,h_x,h_{xx}$:

\begin{align}
\left\{ \begin{array}{ll}
dp_t=&-\Big(\phi_t(\mu_t)p_t+\psi_t(\mu_t)q_t+h_x(x_t,\mu_t)\Big)dt+q_tdB_t\\
p_T=&g_x(x_T) \label{adjoint21}
\end{array} \right.
\end{align}
\begin{align}
\left\{ \begin{array}{ll}
dP_t=&-\Big(2\phi_t(\mu_t)P_t+\psi_t(\mu_t)P_t\psi_t(\mu_t)+2\psi_t(\mu_t)Q_t\\
&+h_{xx}(x_t,\mu_t)\Big)dt+Q_tdB_t\\
P_T=&g_{xx}(x_T) \label{adjoint22}
\end{array} \right.
\end{align}

\medskip
\noindent The $\mathcal{H}$-function associated with $\mu_t$ is defined analogously;
\begin{align}
\mathcal{H}^{(x_t,\mu_t)}(t,x,u)=&H\Big(t,x,u,p_t,q_t-P_t\big(\chi_t(\mu_t)+\psi_t(\mu_t)x_t\big)\Big)\notag\\
&-\frac{1}{2}\Big(\chi_t(u)+\psi_t(u)x\Big)P_t\Big(\chi_t(u)+\psi_t(u)x\Big), \label{h-funct1}
\end{align}
for $(t,x,u)\in [0,T]\times \mathbb{R}\times U$.

\medskip
\noindent Finally, the $\mathcal{H}$-function with  $\nu \in \mathcal{P}(U)$ as the control variable is just (\ref{h-funct1}) integrated with respect to $\nu$:
\begin{align}
\mathcal{H}^{(x_t,\mu_t)}(t,x,\nu)=&H\Big(t,x,\nu,p_t,q_t-P_t\big(\chi_t(\mu_t)+\psi_t(\mu_t)x_t\big)\Big)\notag\\
&-\frac{1}{2}\Big(\chi_t(\nu)+\psi_t(\nu)x\Big)P_t\Big(\chi_t(\nu)+\psi_t(\nu)x\Big), \label{h-funct2}
\end{align}
for $(t,x,\nu) \in [0,T]\times \mathbb{R}\times \mathcal{P}(U)$.

\medskip
\noindent The following result states that the integrated $\mathcal{H}$-function associated with the optimal relaxed control $\hat{\mu}$ is the limit of the integrated $\mathcal{H}$-function associated with the approximating strict control sequence $u^{(k)}$.

\begin{lem}\label{convergenceH}
Let $x_t^{(k)}$ be the state process corresponding to the control sequence $u^{(k)}$ given by the Chattering Lemma, then it holds that
\begin{align*}
\lim_{k\to\infty}\E\int_0^T\mathcal{H}^{(x_t^{(k)},u_t^{(k)})}(t,x_t^{(k)},u_t^{(k)})dt=\E\int_0^T\mathcal{H}^{(\hat{x}_t,\hat{\mu}_t)}(t,\hat{x}_t,\hat{\mu}_t)dt.
\end{align*}
\end{lem}

\medskip
\noindent The main result of this section is the following maximum principle.

\begin{theo}\label{pontryagin} $(Relaxed~ Pontryagin's~ Maximum~ Principle)$\\
If $\hat{\mu}_t$ is an optimal relaxed control with state process $\hat{x}_t$, then for any $t$ outside a null set
\begin{align*}
\mathcal{H}^{(\hat{x}_t,\hat{\mu}_t)}(t,\hat{x}_t,\hat{\mu}_t)&=\sup_{\nu\in\mathcal{P}(U)}\mathcal{H}^{(\hat{x}_t,\hat{\mu}_t)}(t,\hat{x}_t,\nu),~\P\textrm{-a.s.}
\end{align*}
\end{theo}

\medskip\noindent
Using a similar proof as e.g.~Corollaries 4.8 and 4.10 in Bahlali et al.~\cite{bd}, Theorem \ref{pontryagin} is derived from the following integrated maximum principle.

\begin{pro}\label{integratetdmaxprinciple} $(Integrated~ Maximum~ Principle)$ \\
If $\hat{\mu}_t$ is an optimal relaxed control with state process $\hat{x}_t$, then
\begin{align*}
\E\int_0^T\mathcal{H}^{(\hat{x}_t,\hat{\mu}_t)}(t,\hat{x}_t,\hat{\mu}_t)dt=\sup_{u\in U}\E\int_0^T\mathcal{H}^{(\hat{x}_t,\hat{\mu}_t)}(t,\hat{x}_t,u)dt.
\end{align*}
\end{pro}

\noindent\textbf{Proof.}
By the Chattering Lemma there exists an approximating strict control sequence $u_t^{(k)}$ and a corresponding sequence of real numbers $\epsilon^{(k)}\to 0$ such that
\begin{align*}
J(\delta_{u_t^{(k)}})\leq J(\hat{\mu})+\epsilon^{(k)}.
\end{align*}
Sending $k\to\infty$ and using Proposition \ref{nearoptimality} and Lemma \ref{convergenceH} completes the proof.
\eop

\section{An optimal bond portfolio problem} \label{hjmm}

We recall the basic Heath-Jarrow-Morton setup, see Bj\"ork \cite{bjork}. Given a filtered probability space $(\Omega,\mathcal{F},\mathcal{F}_t,\P)$ carrying a $d$-dimensional
$\mathcal{F}_t$-Brownian motion $B_t$, the forward rate $f_t(\tau)$, for each fixed $\tau$, follows an SDE
\begin{align}
df_t(\tau)=\alpha_t(\tau)dt+\tilde{\sigma}_t(\tau)dB_t \label{hjm},
\end{align}
where $\alpha_t(\tau)$ and $\tilde{\sigma}_t(\tau)$ are $\mathbb{R}$- and $\mathbb{R}^d$-valued adapted processes respectively, and $\tau$ denotes time of maturity for the zero
 coupon bond.

\medskip
The bond market induced by the forward rates (\ref{hjm}) is free of arbitrage, in the sense that there exists an equivalent martingale measure, if and only if $\alpha$ can be represented as
\begin{align*}
\alpha_t(\tau)=\tilde{\sigma}_t(\tau)\int_t^\tau\tilde{\sigma}_t(s)ds-\tilde{\sigma}_t(\tau)\Theta_t,
\end{align*}
where $\Theta_t$ is an adapted process such that the Dolean's exponential $\mathcal{E}(\int_0^t \Theta_sdB_s)$ is a $\P$-martingale. The process $\Theta$ is known as the market price
of risk.

\medskip
Consider a market of zero coupon bonds with times \emph{to} maturity in  the interval $U=[0,T^*]$. The forward interest rate under this so-called Musiela parametrization is given by
\begin{align*}
r_t(u)= f_t(t+u), ~ u\in U,
\end{align*}
e.g.~$r_t^0=r_t(0)$ denotes the short rate at time $t$. The re-parametrization yields
\begin{align}
dr_t(u)=\left(\frac{\partial}{\partial u}r_t(u)+\sigma_t(u)\int_0^u\sigma_t(x)dx-\sigma_t(u)\Theta_t\right)dt+\sigma_t(u)dB_t, \label{musiela}
\end{align}
where $\sigma_t(u)=\tilde{\sigma}_t(t+u)$. The relation between bond prices and forward interest rates is given by
\begin{align}
p_t(u)=\exp\left(-\int_0^ur_t(x)dx\right), \label{bondprice-intrate}
\end{align}
and thus by applying It\^o's formula on (\ref{bondprice-intrate}) and inserting Eq.~(\ref{musiela}), one can express the dynamics of the bond prices as
\begin{align}
dp_t(u)&=p_t(u)\left(r_t^0-r_t(u)-v_t(u)\Theta_t\right)dt+p_t(u)v_t(u)dB_t, \label{bondprice}
\end{align}
where $v_t(\cdot)=-\int_0^\cdot\sigma_t(x)dx$.

\medskip
Investing in bonds with the price dynamics as above gives the opportunity to, at any time, choose among a continuum of assets. Namely one for each maturity
 $u\in U$. This gives rise to the problem of how to define a portfolio. A reasonable choice is to consider measure-valued portfolios, as is done in Bj\"ork et al. \cite{bkr}. Using
 measure-valued portfolios also ensures the existence of a locally risk free bank account, $b_t=\exp\left(\int_0^tr_s^0ds\right)$, since this investment is equivalent to a so called roll-over
  strategy, see e.g. Bj\"ork \cite{bjork}. This
 strategy is performed by continuously reinvesting the entire portfolio value in the just maturing bond, and over an arbitrary time interval uses an
  infinite number of assets.

  \medskip
  We define a portfolio as a measure-valued process $\rho_t(du), u\in U$. Intuitively $\rho_t(du)$ is the "number" of bonds in our portfolio at time $t$, with time to maturity
  in the infinitesimal time interval $[u,u+du]$. We denote by $x_t$ the value of the portfolio at time $t$, i.e.
  \begin{align*}
  x_t=\int_Up_t(u)\rho_t(du).
  \end{align*}
 Further, the portfolio is self financing, i.e. the increments of the portfolio value are due to price changes only.
  Referring to Ekeland and Taflin \cite{ekelandtaflin}, we may formally express this as
\begin{align*}
dx_t&=\int_U\rho_t(du)\left(dp_t(u)-\frac{\partial}{\partial u}p_t(u)dt\right),
\end{align*}
where the last term appears because we use the Musiela parametrization. In this setting the portfolio value changes due to both price changes as well as to changes in time to
maturity.
Using Eq.~(\ref{bondprice}) and noting that $\frac{\partial}{\partial u}p_t(u)=-p_t(u)r_t(u)$, the above relation is interpreted as follows (see e.g. Bj\"ork et al. \cite{bkr}).
\begin{align*}
 x_t=&x_0+\int_0^t\int_Up_s(u)\left(r_s^0-r_s(u)-v_s(u)\Theta_s\right)\rho_s(du)ds\\
&+\int_0^t\int_Up_s(u)v_s(u)\rho_s(du)dB_s+\int_0^t\int_Up_s(u)r_s(u)\rho_s(du)ds\\
=&x_0+\int_0^t\int_Up_s(u)\left(r_s^0-v_s(u)\Theta_s\right)\rho_s(du)ds
+\int_0^t\int_Up_s(u)v_s(u)\rho_s(du)dB_s,
\end{align*}
where $x_0$ is the initial capital.
  Considering only portfolios with positive holdings (i.e. no short positions), we may write
  \begin{align}
  \rho_t(du)=\frac{x_t}{p_t(u)}\mu_t(du), \label{relativepf}
  \end{align}
 where $\mu_t\in\boldsymbol{M}(\Omega)$, i.e.~it takes values in the set of probability measures on $U$. This so-called relative portfolio is the proportion of the
  portfolio invested in bonds with time to maturity in $[u,u+du]$.

\medskip
   Now, inserting the expression (\ref{relativepf}) into the portfolio dynamics above yields
\begin{align}
x_t=&x_0+\int_0^tx_s\int_U\left(r_s^0-v_s(u)\Theta_s\right)\mu_s(du)ds+\int_0^tx_s\int_Uv_s(u)\mu_s(du)dB_s. \label{pfdynamics}
\end{align}

\medskip
\noindent The aim is to control this self-financing portfolio, via $\mu_t$, in an optimal way.
Assuming that our goal is to minimize the cost functional
\begin{align*}
J(\mu)=\E\left(\int_0^T\int_Uh(t,x_t,u)\mu_t(du)dt+g(x_T)\right),
\end{align*}
we get an optimal control problem on the form (\ref{linearcoeff}),(\ref{sde1}),(\ref{costfunctional}), with $\upsilon_t=\chi_t\equiv 0$,
$\phi_t(u)=\big(r_t^0-v_t(u)\Theta_t\big)$ and $\psi_t(u)=v_t(u)$.\\

\begin{ex}\label{example1}

Consider a passive investor who invests the initial capital into $N$ number of bonds with times \textbf{of} maturity $T_1,\ldots,T_N$, at time $0$ and does nothing thereafter, i.e.
\begin{align*}
x_t=\sum_{i=1}^N\tilde{p}_{t}(T_i),
\end{align*}
where $\tilde{p}_t(T)$ denotes the price of a bond with time of maturity $T$. This corresponds to a portfolio
consisting of bonds which, at time $t$, have times to maturity $T_1-t,\ldots,T_N-t$. Thus, our relative portfolio is
\begin{align*}
\mu_t(du)=\sum_{i=1}^N\delta_{(T_i-t)}(du),
\end{align*}
and \textnormal{Eq.~(\ref{pfdynamics})} becomes
\begin{align*}
x_t=&\sum_{i=1}^Np_{0}(T_i)+\int_0^tx_s\sum_{i=1}^N\left(r_s^0-v_s(T_i-s)\Theta_s\right)ds
+\int_0^tx_s\sum_{i=1}^Nv_s(T_i-s)dB_s.
\end{align*}
Note that by the Musiela parametrization we are in a moving time frame and therefore, although the investor is passive, the control measure changes continuously in $t$.
\end{ex}

\subsection{Mean variance portfolio selection}
In this last section we derive the adjoint equations and $\mathcal{H}$-function for a specific example. The cost functional corresponds to a mean variance portfolio selection problem and two different choices of
interest rate processes are considered. In principle, necessary conditions for a portfolio to be optimal can be found by maximizing the
$\mathcal{H}$-function with respect to a measure on $U$. Unfortunately, the BSDEs for $(p,q)$ and $(P,Q)$ are quite involved and it seems difficult to find explicit solutions.\\

\noindent Assume the following cost functional:
\begin{align}
J(\mu)=\E\left(\frac{1}{2}(x_T-\kappa)^2\right), \label{meanvar}
\end{align}
with given constant $\kappa$. Minimizing $J$ (for a certain $\kappa$) is equivalent to a mean variance portfolio selection problem.
 Assume that $\hat{\mu}_t$ with corresponding portfolio value $\hat{x}_t$ is optimal. Using the relaxed maximum principle we may write down the necessary conditions for
 $\hat{\mu}_t$ and $\hat{x}_t$. The adjoint equations becomes
\begin{align*}
\left\{ \begin{array}{ll}
dp_t=&-\Big((r_t^0-v_t(\hat{\mu}_t)\Theta_t)p_t+v_t(\hat{\mu}_t)q_t\Big)dt+q_tdB_t,\\
p_T=&\hat{x}_T-\kappa,
\end{array} \right.
\end{align*}
and
\begin{align*}
\left\{ \begin{array}{ll}
dP_t=&-\Big(\big(2(r_t^0-v_t(\hat{\mu}_t)\Theta_t)+(v_t(\hat{\mu}_t))^2\big)P_t+2v_t(\hat{\mu}_t)Q_t\Big)dt+Q_tdB_t,\\
P_T=&1.
\end{array} \right.
\end{align*}
The corresponding $\mathcal{H}$-function is
\begin{align}
\mathcal{H}^{(\hat{x}_t,\hat{\mu}_t)}(t,\hat{x}_t,\nu)=&-p_t\Big(r_t^0-v_t(\nu)\Theta_t\Big)\hat{x}_t-\Big(q_t-P_tv_t(\hat{\mu}_t)\hat{x}_t\Big)v_t(\nu)\hat{x}_t\notag\\
&-\frac{1}{2}P_tv_t^2(\nu)\hat{x}_t^2, \label{hfunct}
\end{align}
where again, $v_t(\nu)=\int_Uv_t(u)\nu(du)$ and $v_t^2(\nu)=\int_Uv_t^2(u)\nu(du)$.
\subsubsection{Ho-Lee}
Choosing the volatility process to be constant,
\begin{align*}
\sigma_t(u)=\sigma,
\end{align*}
and consequently
\begin{align*}
v_t(u)=-\sigma u,
\end{align*}
the short rate $r_t^0$ is a Gaussian process. Under a obvious integrability assumption on $\Theta_t$ we then have that (A.1)-(A.3) are fulfilled. Thus by the relaxed maximum principle, a necessary condition for a portfolio $\hat{\mu}_t$ to
minimize the cost functional (\ref{meanvar}), is that it maximizes Eq.~(\ref{hfunct}) with
\begin{align*}
v_t(\nu)=-\sigma\int_Uu\nu(du).
\end{align*}

\subsubsection{Hull-White}
Another choice of volatility process that induces a mean-reverting Gaussian short rate is
\begin{align*}
\sigma_t(u)=\sigma e^{-cu},
\end{align*}
with constants $\sigma$ and $c$. Thus, a necessary condition for optimality of a portfolio $\hat{\mu}_t$ is that it maximizes Eq.~(\ref{hfunct}) with
\begin{align*}
v_t(\nu)=\frac{\sigma}{c}\int_U(e^{-cu}-1)\nu(du).
\end{align*}


\section{Proofs and technical results}
Throughout this section, we denote by $K>0$ a generic constant that may vary from line to line.\\

\medskip
\noindent\textbf{Proof of Lemma \ref{contmeasure}.}
The second assertion follows easily from the first by  using the Lipschitz property
of $g$. For the first assertion, 
define (pointwise in $\Omega$) for $k=1,2,\ldots$ and $A$ a Borel set; 
\begin{align*}
\mu_t^{(k)}(A)=\left\{ \begin{array}{ll} 
\frac{1}{t}\int_0^t\mu_s(A)ds, & t\in(0,2^{-k})\\
& \\
2^k\int_{t-2^{-k}}^t\mu_s(A)ds, & t\in[2^{-k},T).
\end{array} \right.
\end{align*}
Obviously, for each $k$, $\mu_t^{(k)}(A)$ is continuous for all Borel sets $A$.
Moreover, one can show (cf. Ma and Yong \cite{mayong}) that for any $f\in
C(U,\mathbb{R}^n)$ the following holds for $(t,\omega)$ outside a null set
\begin{align*}
f(\mu_t^{(k)})\to f(\mu_t), ~ \textrm{as} ~ k\to\infty.
\end{align*} 
In particular,
\begin{align}
\varphi_t(\mu_t^{(k)})\to \varphi_t(\mu_t), \label{convas}
\end{align}
$\P\textrm{-a.s.}$ for $t$ outside a null set, where $\varphi$ stands for one of the
processes $\upsilon, \phi, \chi, \psi$. Denote $y_t=x_t^{(k)}-x_t$, then $y_t$ can
be expressed as 
\begin{align*}
y_t=&\int_0^t\Big(\nu_s(\mu_s^{(k)})-\nu_s(\mu_s)+\big(\phi_s(\mu_s^{(k)})-\phi_s(\mu_s)\big)x_s^{(k)}\Big)ds+\int_0^t\phi_s(\mu_s)y_sds\\
+&\int_0^t\Big(\chi_s(\mu_s^{(k)})-\chi_s(\mu_s)+\big(\psi_s(\mu_s^{(k)})-\psi_s(\mu_s)\big)x_s^{(k)}\Big)dB_s+\int_0^t\psi_s(\mu_s)y_sdB_s.
\end{align*}
Let
\begin{align*}
z_t=1-\int_0^t\phi_s(\mu_s)z_sds,
\end{align*}
and apply Ito's formula on $z_ty_t$ to get
\begin{align*}
z_ty_t=&\int_0^t\Big(\nu_s(\mu_s^{(k)})-\nu_s(\mu_s)+\big(\phi_s(\mu_s^{(k)})-\phi_s(\mu_s)\big)x_s^{(k)}\Big)z_sds\\
&+\int_0^t\Big(\chi_s(\mu_s^{(k)})-\chi_s(\mu_s)+\big(\psi_s(\mu_s^{(k)})-\psi_s(\mu_s)\big)x_s^{(k)}\Big)z_sdB_s\\
&+\int_0^t\psi_s(\mu_s)z_sy_sdB_s.
\end{align*}
Noting that $\psi$ is bounded, we apply the Burkholder-Davis-Gundy, Gronwall and H\"older inequalities to get that for any $p\geq1$,
\begin{align*}
&\E|zy|_T^{*,2p}\\
\leq &K\bigg\{\E\int_0^T\Big(|\nu_s(\mu_s^{(k)})-\nu_s(\mu_s)|^{2p}+|\phi_s(\mu_s^{(k)})-\phi_s(\mu_s)|^{2p}x_t^{2p}\Big)dt\\
&+\E\int_0^T\Big(|\chi_s(\mu_s^{(k)})-\chi_s(\mu_s)|^{2p}+|\psi_s(\mu_s^{(k)})-\psi_s(\mu_s)|^{2p}x_t^{2p}\Big)dt\bigg\}\\
\leq&K\bigg\{\E\int_0^T|\nu_s(\mu_s^{(k)})-\nu_s(\mu_s)|^{2p}dt+\Big(\E\int_0^T|\phi_s(\mu_s^{(k)})-\phi_s(\mu_s)|^{4p}dt\Big)^{1/2}\\
&+\E\int_0^T|\chi_s(\mu_s^{(k)})-\chi_s(\mu_s)|^{2p}dt+\Big(\E\int_0^T|\psi_s(\mu_s^{(k)})-\psi_s(\mu_s)|^{4p}dt\Big)^{1/2}\bigg\},
\end{align*}
where all the terms on the right hand side converge to $0$ as $k\to\infty$ by the
Dominated Convergence Theorem. We note that by (A.3) we have $\E|z|_T^{*,p}<\infty$, $\E|z^{-1}|_T^{*,p}<\infty$, for any $p\geq 1$. Thus, using the H\"older inequality
we conclude that
\begin{align*}
\E|x^{(k)}-x|_T^{*,2}=\E|z^{-1}zy|_T^{*,2}\leq K\Big(\E|zy|_T^{*,4}\Big)^{1/2}\to0,
\end{align*}
as $k\to\infty$.
\eop

\medskip
 \noindent\textbf{Proof of Lemma \ref{tight}.} It suffices to check that the marginals are tight. $\boldsymbol{B}_t^{(k)}$ is tight since the processes induce the same measure for every $k$. Further,
$\lambda_t^{(k)}$ is tight because $\mathcal{L}([0,T]\times U)$ is compact. Finally by (\ref{kolmogorovcondition}), there exists a constant $K>0$ such that
\begin{align*}
\E^{(k)}(x_t^{(k)}-x_s^{(k)})^4\leq K|s-t|^2,
\end{align*}
for all $t,s \in [0,T]$, for all $k$, where $\E^{(k)}$ is the expectation under $\P^{(k)}$.
 Hence the Kolmogorov condition is fulfilled (see e.g. Yong and Zhou \cite{yongzhou}, Theorem 2.14.) and $x_t^{(k)}$ is tight.
\eop

\medskip
\noindent\textbf{Proof of Theorem \ref{existence}.}
Define the set $\mathcal{K}=\Big\{\big(b(s,x_s,\mu_s),\sigma(s,x_s,\mu_s)\big):\mu_s\in\boldsymbol{M}(\Omega)\Big\}$. Then $\mathcal{K}$ is a convex set in
$L^2([0,T]\times\Omega)\times L^2([0,T]\times\Omega)$ and by Mazur's Theorem the weak closure of $\mathcal{K}$ equals its strong closure. Thus, by Lemma \ref{weakconv}, for
each integer $l>0$ and $\epsilon>0$, there exists a finite set of numbers $\{\alpha_1,\ldots,\alpha_{N(l,\epsilon)}\}$ satisfying $\alpha_i\geq 0$ and $\sum_i\alpha_i=1$, such
that \begin{align}
&\E\int_0^T\big|\sum_{i=1}^{N(l,\epsilon)}\alpha_ib(t,x_t,\mu_t^{(l+i)})-b(t,x_t,\mu_t)\big|^2dt\\
&+\E\int_0^T\big|\sum_{i=1}^{N(l,\epsilon)}\alpha_i\sigma(t,x_t,\mu_t^{(l+i)})-\sigma(t,x_t,\mu_t)\big|^2dt<\epsilon.
\label{L2ineq}
 \end{align}
By uniform integrability, for any $\epsilon >0$, there exists an integer $N_0(\epsilon)>0$, such that
\begin{align*}
\E\Big(&\big|x^{(k)}-x\big|_T^{*,4}+\big|B^{(k)}-B\big|_T^{*,4}+\big|\upsilon^{(k)}-\upsilon\big|_T^{*,4}\\
&+\big|\phi^{(k)}-\phi\big|_T^{*,4}+\big|\chi^{(k)}-\chi\big|_T^{*,4}+\big|\psi^{(k)}-\psi\big|_T^{*,4}\Big)<\epsilon^2,
\end{align*}
 for all $k>N_0$ and $u\in U$.
Fix such an $\epsilon$ and let $\bar{N}=N(N_0,\epsilon)$ and $\{\alpha_1,\ldots,\alpha_{\bar{N}}\}$ be such that $\alpha_i\geq0$; $\sum_i\alpha_i=1$ and (\ref{L2ineq}) holds.
Denote by $K$ a generic constant that may vary from line to line.
Define for each $k=1,2,\ldots$,
\begin{align*}
b^{(k)}(t,x_t,\mu_t)=\upsilon_t^{(k)}(\mu_t)+\phi_t^{(k)}(\mu_t)x_t,\\
\sigma^{(k)}(t,x_t,\mu_t)=\chi_t^{(k)}(\mu_t)+\psi_t^{(k)}(\mu_t)x_t,
\end{align*}
and for each $i=1,\ldots,\bar{N}$
\begin{align*}
\Delta^i(B)_t=&\int_0^t\sigma^{(N_0+i)}(s,x_s^{(N_0+i)},\mu_s^{(N_0+i)})dB_s^{(N_0+i)}\\
&-\int_0^t\sigma^{(N_0+i)}(s,x_s^{(N_0+i)},\mu_s^{(N_0+i)})dB_s.
\end{align*}
Then it is readily seen that
\begin{align*}
\E\big|\sum_i^{\bar{N}}\alpha_i\Delta^i(B)\big|_T^{*,2}<K\epsilon .
\end{align*}
Similarly, let
\begin{align*}
\Delta^i(x)^b_t &= \int_0^tb^{(N_0+i)}(s,x_s^{(N_0+i)},\mu_s^{(N_0+i)})ds
-\int_0^tb^{(N_0+i)}(s,x_s,\mu_s^{(N_0+i)})ds,\\
\Delta^i(x)^\sigma_t &= \int_0^t\sigma^{(N_0+i)}(s,x_s^{(N_0+i)},\mu_s^{(N_0+i)})dB_s
-\int_0^t\sigma^{(N_0+i)}(s,x_s,\mu_s^{(N_0+i)})dB_s,\\
\Delta^i(b)_t&=\int_0^tb^{(N_0+i)}(s,x_s,\mu_s^{(N_0+i)})ds
-\int_0^tb(s,x_s,\mu_s^{(N_0+i)})ds,\\
\Delta^i(\sigma)_t&=\int_0^t\sigma^{(N_0+i)}(s,x_s,\mu_s^{(N_0+i)})dB_s
-\int_0^t\sigma(s,x_s,\mu_s^{(N_0+i)})dB_s.
\end{align*}
and conclude that
\begin{align*}
\E&\Big(\big|\sum_i^{\bar{N}}\alpha_i\Delta^i(x)^b\big|_T^{*,2}+\big|\sum_i^{\bar{N}}\alpha_i\Delta^i(x)^\sigma\big|_T^{*,2}
+\big|\sum_i^{\bar{N}}\alpha_i\Delta^i(b)\big|_T^{*,2}+\big|\sum_i^{\bar{N}}\alpha_i\Delta^i(\sigma)\big|_T^{*,2}\Big)\\
&<K\epsilon,
\end{align*}
where we have used the Burkholder-Davis-Gundy inequality for the martingale terms. \\

Note that for each $k$, $x_t^{(k)}$ satisfies
\begin{align*}
x_t^{(k)}=x_0+\int_0^tb^{(k)}(s,x_s^{(k)},\mu_s^{(k)})ds+\int_0^t\sigma^{(k)}(s,x^{(k)}_s,\mu_s^{(k)})dB^{(k)}_s,
\end{align*}  
and thus 
\begin{align*}
\sum_i^{\bar{N}}\alpha_ix^{(N_0+i)}_t&=x+\int_0^tb(s,x_s,\mu)ds+\int_0^t\sigma(s,x_s,\mu)dB_s\\
+&\sum_i^{\bar{N}}\Delta^i(B)_t+\sum_i^{\bar{N}}\Delta^i(x)^b_t+\sum_i^{\bar{N}}\Delta^i(x)^\sigma_t
+\sum_i^{\bar{N}}\Delta^i(b)_t+\sum_i^{\bar{N}}\Delta^i(\sigma)_t\\
+&\int_0^t\bigg(\sum_i^{\bar{N}}\alpha_ib(s,x_s,\mu_s^{(N_0+i)})-b(s,x_s,\mu_s)\bigg)ds\\
+&\int_0^t\bigg(\sum_i^{\bar{N}}\alpha_i\sigma(s,x_s,\mu_s^{(N_0+i)})-\sigma(s,x_s,\mu_s)\bigg)dB_s. 
\end{align*} 
Since $f(x)=x^2$ is
convex, it is easy to check that \begin{align*} \E\bigg(\big|\sum_i^{\bar{N}}\alpha_ix^{(N_0+i)}-x\big|_T^{*,2}\bigg)\leq \epsilon. 
\end{align*} Combining this with
the previous inequalities and using the Burkholder-Davis-Gundy inequality, yields
\begin{align*} 
&\E\Big(\big|x_t-x_0-\int_0^tb(s,x_s,\mu)ds-\int_0^t\sigma(s,x_s,\mu)dB_s\big|_T^{*,2}\Big)\\
\leq&
K\bigg\{\epsilon+\E\int_0^T\big|\sum_{i=1}^{\bar{N}}\alpha_ib(s,x_s,\mu_s^{(N_0+i)})-b(s,x_s,\mu)\big|^2dt\\
&+\E\int_0^T\big|\sum_{i=1}^{\bar{N}}\alpha_i\sigma(s,x_s,\mu_s^{(N_0+i)})
-\sigma(s,x_s,\mu)\big|^2dt\bigg\}\\
<&K\epsilon.
\end{align*}
This completes the proof. \eop

\medskip
\noindent\textbf{Proof of Lemma \ref{chattering}.}
The construction of the approximating sequence $u_t^{(k)}$ is done as in Ma and Yong \cite{mayong}, Theorem 3.6. Then by tightness and Skorohod's Selection Theorem
$\delta_{\hat{u}_t^{(k)}}(du)dt\to\hat{\mu}_t(du)dt$ $\hat{\P}$-a.s. in
$\mathcal{L}([0,T]\times U)$ and (\ref{chattering1}) follows in the same way as in Theorem \ref{existence}.
\eop

\medskip
\noindent\textbf{Proof of Proposition \ref{varineq}.}
Let $x^\theta_t$ denote the state process corresponding to $u^\theta$. We proceed as in Peng \cite{peng} to introduce the first- and second order variational equations (noting that $b_x=\phi, \sigma_x=\psi, b_{xx}=\sigma_{xx}=0$ in our case):
\begin{align*}
x^{(1)}_t=&\int_0^t\Big(\phi_s(u_s)x^{(1)}_s+\nu_s(u^\theta_s)+\phi_s(u^\theta_s)x_s-\nu_s(u_s)-\phi_s(u_s)x_s\Big)ds\\
&+\int_0^t\Big(\psi_s(u_s)x_s^{(1)}+\chi_s(u^\theta_s)+\psi_s(u^\theta_s)x_s-\chi_s(u_s)-\psi_s(u_s)x_s\Big)dB_s,
\end{align*}
\begin{align*}
x^{(2)}_t=&\int_0^t\Big(\phi_s(u_s)x^{(2)}_s+\big(\phi_s(u^\theta_s)-\phi_s(u_s)\big)x^{(1)}_s\Big)ds\\
&+\int_0^t\Big(\psi_s(u_s)x^{(2)}_s+\big(\psi_s(u^\theta_s)-\psi_s(u_s)\big)x_s^{(1)}\Big)dB_s.
\end{align*}
Then we have the following estimate.
\begin{align}
\E|x^\theta-x-x^{(1)}-x^{(2)}|_T^{*,2}\leq K\theta^2.\label{varineq1}
\end{align}
To prove (\ref{varineq1}), note that for any $p\geq 1$ it holds that
\begin{align}
\E|x^{(1)}|_T^{*,2p}\leq K\theta^p \label{varineq2}
\end{align}
and
\begin{align}
\E|x^{(2)}|_T^{*,2p}\leq K\theta^{2p}. \label{varineq3}
\end{align}
As in Peng \cite{peng}, we can write
\begin{align*}
x_t+x_t^{(1)}+x_t^{(2)}=&x_0+\int_0^t\Big(\nu_s(u^\theta_s)+\phi_s(u^\theta_s)(x_s+x_s^{(1)}+x_s^{(2)})\Big)ds\\
&+\int_0^t\Big(\chi_s(u^\theta_s)+\psi_s(u^\theta_s)(x_s+x_s^{(1)}+x_s^{(2)})\Big)dB_s\\
&-\int_0^tG^\theta_sds-\int_0^t\Lambda^\theta_sdB_s.
\end{align*}
where
\begin{align*}
G^\theta_s=\Big(\phi_s(u^\theta_s)-\phi_s(u_s)\Big)x_s^{(2)},
\end{align*}
\begin{align*}
\Lambda^\theta_s=\Big(\psi_s(u^\theta_s)-\psi_s(u_s)\Big)x_s^{(2)}.
\end{align*}
We have for $G^\theta$ and $\Lambda^\theta$ that
\begin{align}
\E\Big|\int_0^\cdot G^\theta_sds\Big|_T^{*,2}+\E\Big|\int_0^\cdot \Lambda^\theta_sdB_s\Big|_T^{*,2}\leq K\theta^2. \label{varineq8}
\end{align}
Thus,
\begin{align*}
x^\theta_t-x_t-x_t^{(1)}-x_t^{(2)}&=\int_0^t\phi_s(u^\theta_s)(x^\theta_s-x_s-x_s^{(1)}-x_s^{(2)})ds\\
&+\int_0^t\psi_s(u^\theta_s)(x^\theta_s-x_s-x_s^{(1)}-x_s^{(2)})dB_s\\
&+\int_0^tG^\theta_sds+\int_0^t\Lambda_s^\theta dB_s,
\end{align*}
which together with (\ref{varineq2}), (\ref{varineq3}) and (\ref{varineq8}) leads to (\ref{varineq1}). With this result, a Taylor expansion of the cost functional as in Peng \cite{peng} gives us
\begin{align}
&J(\delta_{u^\theta})-J(\delta_{u})\notag\\
=&\E\int_0^T\Big(h_x(x_s,u_s)(x^{(1)}_s+x^{(2)}_s)+\frac{1}{2}h_{xx}(x_s,u_s)x_s^{(1)}x_s^{(1)}\Big)ds\notag\\
+&\E\int_0^T\Big(h(x_s,u^\theta_s)-h(x_s,u_s)\Big)ds 
+\E\Big(g_x(x_T)(x_T^{(1)}+x_T^{(2)})\Big)\label{varineq5} \\
+&\frac{1}{2}\E\Big(g_{xx}(x_T)x_T^{(1)}x_T^{(1)}\Big)+o(\theta).\notag 
\end{align}
The next step is to express the right hand side in terms of the first- and second order adjoint processes. We start by deriving the former. Define
\begin{align*}
\Phi_t=1+\int_0^t\phi_s(u_s)\Phi_sds+\int_0^t\psi_s(u_s)\Phi_sdB_s.
\end{align*}
By Ito's formula $\Phi^{-1}$ is given by
\begin{align*}
\Phi_t^{-1}=1+\int_0^t\big(\psi_s(u_s)\psi_s(u_s)-\phi_s(u_s)\big)\Phi^{-1}_sds-\int_0^t\psi_s(u_s)\Phi^{-1}_sdB_s.
\end{align*}
By a simple manipulation we deduce the moment property 
\begin{align}
\E|\Phi|_T^{*,p}+\E|\Phi^{-1}|_T^{*,p}<\infty, \label{varineq4}
\end{align}
for any $p\geq1$. Next, we introduce
\begin{align*}
X^{(1)}&=\Phi_Tg_x(x_T)+\int_0^T\Phi_th_x(x_t,u_t)dt,\\
y^{(1)}_t&=\E\big(X^{(1)}|\mathcal{F}_t\big)-\int_0^t\Phi_sh_x(x_s,u_s)ds.
\end{align*}
Since $g_x$ and $h_x$ are bounded we can use (\ref{varineq4}) to deduce that
\begin{align*}
\E\big|X^{(1)}\big|^p<\infty,
\end{align*}
for any $p\geq1$. Thus by the Martingale Representation Theorem there exists an $\mathcal{F}_t$-adapted process $H_t$ with the property that, for any $p\geq1$,
\begin{align*}
\E\int_0^T|H_t|^pdt<\infty,
\end{align*}
and such that
\begin{align*}
y^{(1)}_t=\E(X^{(1)})+\int_0^tH_sdB_s-\int_0^t\Phi_sh_x(x_s,u_s)ds.
\end{align*}
We may now define our first order adjoint processes $(p,q)$ as
\begin{align*}
p_t&=\Phi_t^{-1}y^{(1)}_t,\\
q_t&=\Phi_t^{-1}H_t-\psi_t(u_t)p_t,
\end{align*}
noting that $p_t$ and $q_t$ are $\mathbb{R}$- resp. $\mathbb{R}^d$-valued $\mathcal{F}_t$-adapted processes satisfying
\begin{align*}
\E|p|_T^{*,p}+\E\int_0^T|q|^pdt<\infty,
\end{align*}
for any $p\geq1$. Applying Ito's formula on $p_t=\Phi_t^{-1}y^{(1)}_t$ yields
\begin{align}
dp_t=-\Big(h_x(x_t,u_t)+p_t\phi_t(u_t)+q_t\psi_t(u_t)\Big)dt+q_tdB_t.\label{varineq10}
\end{align}
Using (\ref{varineq10}), and once again by using Ito's formula we can derive
\begin{align*}
\E\Big(p_Tx_T^{(1)}\Big)=&\E\int_0^T\Big(p_t\big(\nu_t(u_t^\theta)+\phi_t(u_t^\theta)x_t-\nu_t(u_t)-\phi_t(u_t)x_t\big)\\
&+q_t\big(\chi_t(u_t^\theta)+\psi_t(u_t^\theta)x_t-\chi_t(u_t)-\psi_t(u_t)x_t\big)-h_x(x_t,u_t)x_t^{(1)}\Big)dt.
\end{align*}
Similarly,
\begin{align}
&\E\Big(p_Tx_T^{(2)}\Big)\notag\\
=&\E\int_0^T\Big(p_t\big(\phi_t(u_t^\theta)-\phi_t(u_t)\big)x_t^{(1)}+q_t\big(\psi_t(u_t^\theta)-\psi_t(u_t)\big)x_t^{(1)}-h_x(x_t,u_t)x_t^{(2)}\Big)dt. \label{varineq6}
\end{align}
Now, noting that $p_T=\Phi_T^{-1}\Phi_Tg_x(x_T)=g_x(x_T)$ and that the first two terms on the right hand side of (\ref{varineq6}) is of order $\theta^{3/2}$, we may combine the two equalities above to
\begin{align*}
\E\Big(g_x(x_T)(x_T^{(1)}+x_T^{(2)})\Big)=&\E\int_0^Tp_t\Big(\nu_t(u_t^\theta)+\phi_t(u_t^\theta)x_t-\nu_t(u_t)-\phi_t(u_t)x_t\Big)dt\\
&+\E\int_0^Tq_t\Big(\chi_t(u_t^\theta)+\psi_t(u_t^\theta)x_t-\chi_t(u_t)-\psi_t(u_t)x_t\Big)dt\\
&-\E\int_0^T\Big(h_x(x_t,u_t)(x_t^{(1)}+x_t^{(2)})\Big)dt+o(\theta).
\end{align*}
Thus, we may write the right hand side of (\ref{varineq5}) as
\begin{align*}
&\E\int_0^T\Big(H(t,x_t,u_t,p_t,q_t)-H(t,x_t,u^\theta_t,p_t,q_t)\Big)dt+\frac{1}{2}\E\int_0^Th_{xx}(x_t,u_t)x^{(1)}_tx^{(1)}_tdt\\
&+\frac{1}{2}\E\Big(g_{xx}(x_T)x^{(1)}_Tx^{(1)}_T\Big)+o(\theta).
\end{align*}
It remains to replace the second order terms with the second order adjoint processes. Define
\begin{align*}
\Psi_t=1+\int_0^t\big(2\phi_s(u_s)+\psi_s(u_s)\psi_s(u_s)\big)\Psi_sds+\int_0^t2\psi_s(u_s)\Psi_sdB_s.
\end{align*}
Thus, $\Psi^{-1}$ is given by
\begin{align*}
\Psi_t^{-1}=1+\int_0^t\big(3\psi_s(u_s)\psi_s(u_s)-2\phi_s(u_s)\big)\Psi^{-1}_sds-\int_0^t2\psi_s(u_s)\Psi^{-1}_sdB_s,
\end{align*}
and as above the moment property
\begin{align*}
\E|\Phi|_T^{*,p}+\E|\Phi^{-1}|_T^{*,p}<\infty,
\end{align*}
holds for any $p\geq1$. Next, we introduce
\begin{align*}
X^{(2)}&=\Psi_Tg_{xx}(x_T)+\int_0^T\Psi_th_{xx}(x_t,u_t)dt,\\
y^{(2)}_t&=\E\big(X^{(2)}|\mathcal{F}_t\big)-\int_0^t\Psi_sh_{xx}(x_s,u_s)ds,
\end{align*}
where $g_{xx}$ and $h_{xx}$ are bounded so that
\begin{align*}
\E\big|X^{(2)}\big|^p<\infty,
\end{align*}
for any $p\geq1$. Again, by the Martingale Representation Theorem there exists an $\mathcal{F}_t$-adapted process $K_t$ with the property that, for any $p\geq1$,
\begin{align*}
\E\int_0^T|K_t|^pdt<\infty,
\end{align*}
and such that
\begin{align*}
y^{(2)}_t=\E(X^{(2)})+\int_0^tK_sdB_s-\int_0^t\Psi_sh_{xx}(x_s,u_s)ds.
\end{align*}
The second order adjoint processes $(P,Q)$ is defined as
\begin{align*}
P_t&=\Psi_t^{-1}y^{(2)}_t,\\
Q_t&=\Psi_t^{-1}K_t-2\psi_t(u_t)P_t,
\end{align*}
noting that $P_t$ and $Q_t$ are $\mathbb{R}$- resp. $\mathbb{R}^d$-valued $\mathcal{F}_t$-adapted processes satisfying
\begin{align*}
\E|P|_T^{*,p}+\E\int_0^T|Q|^pdt<\infty,
\end{align*}
for any $p\geq1$.  Applying Ito's formula on $P_t=\Psi_t^{-1}y^{(2)}_t$ yields
\begin{align*}
dP_t=-\Big(h_{xx}(x_t,u_t)+2P_t\phi_t(u_t)+P_t\psi_t(u_t)\psi_t(u_t)+2Q_t\psi_t(u_t)\Big)dt+Q_tdB_t.
\end{align*}
By another application of Ito's formula we deduce
\begin{align}
&\E\Big(P_Tx_T^{(1)}x_T^{(1)}\Big)\notag\\
=&\E\int_0^T\Big(2P_t\big(b(x_t,u_t^\theta)-b(x_t,u_t)+\psi_t(u_t)(\sigma(x_t,u_t^\theta)-\sigma((x_t,u_t)\big)x_t^{(1)}\notag\\
&+2Q_t\big(\sigma(x_t,u_t^\theta)-\sigma((x_t,u_t)\big)x_t^{(1)}-h_{xx}(x_t,u_t)x_t^{(1)}x_t^{(1)}\label{varineq7}\\
&+\big(\sigma(x_t,u_t^\theta)-\sigma(x_t,u_t)\big)P_t\big(\sigma(x_t,u_t^\theta)-\sigma(x_t,u_t)\big)\Big)dt.\notag
\end{align}
Noting that $P_T=\Phi_T^{-1}\Phi_Tg_{xx}(x_T)=g_{xx}(x_T)$ and that the first two terms on the right hand side of (\ref{varineq7}) is of order $\theta^{3/2}$, we get
\begin{align*}
\E\Big(g_{xx}(x_T)(x_T^{(1)}x_T^{(1)})\Big)
=&\E\int_0^T\big(\sigma(x_t,u_t^\theta)-\sigma(x_t,u_t)\big)P_t\big(\sigma(x_t,u_t^\theta)-\sigma(x_t,u_t)\big)dt\\
&-\E\int_0^T\Big(h_{xx}(x_t,u_t)x_t^{(1)}x_t^{(1)}\Big)dt+o(\theta).
\end{align*}
Hence we can remove the second order terms on the right hand side of (\ref{varineq5}) to get
\begin{align*}
&\E\int_0^T\Big(H(t,x_t,u_t,p_t,q_t)-H(t,x_t,u^\theta_t,p_t,q_t)\Big)dt\\
&+\frac{1}{2}\E\int_0^T\big(\sigma(x_t,u_t^\theta)-\sigma(x_t,u_t)\big)P_t\big(\sigma(x_t,u_t^\theta)-\sigma(x_t,u_t)\big)dt+o(\theta)\\
&=\E\int_0^T\Big(\mathcal{H}^{(x_t,u_t)}(t,x_t,u_t)-\mathcal{H}^{(x_t,u_t)}(t,x_t,u^\theta_t)\Big)dt+o(\theta),
\end{align*}
which completes the proof.
\eop

Before we can prove Proposition \ref{nearoptimality}, we need some preliminary results.
We start by defining a metric on $\mathcal{U}$:
\begin{align}
d(u,u')=\P\otimes dt\Big\{(\omega,t)\in\Omega\times [0,T];u(\omega,t)\neq u'(\omega,t)\Big\},
\end{align}
where $\P\otimes dt$ is the product measure of $\P$ and the Lebesgue measure. Then we have the following result.

\begin{lem} \label{conttrajectory}
\begin{itemize}
\item[(i)]
$(\mathcal{U},d)$ is a complete metric space.
\item[(ii)]
For any $p\geq 1$ there exists a constant $K$ such that for any $u,u'\in\mathcal{U}$, it holds that
\begin{align}
\E|x-x'|_T^{*,2p}\leq K\big(d(u,u')\big)^{1/4},
\end{align}
where $x_t$ and $x_t'$ are the state processes corresponding to $u$ and $u'$.
\item[(iii)]
For any $u,u'\in\mathcal{U}$ there exists a constant $K$ such that
\begin{align}
|J(\delta_u)-J(\delta_{u'})|\leq K\big(d(u,u')\big)^{1/4}.
\end{align}
\end{itemize}
\end{lem}

\noindent\textbf{Proof.}
$(i)$ can be proved as in Elliott and Kohlmann \cite{elliotkohlmann}, see also Zhou \cite{zhou}.\\
\medskip \noindent As for $(ii)$, denoting $y_t=x_t-x_t'$, we have that $y_t$ satisfies the SDE
\begin{align*}
y_t&=\int_0^t\Big(\nu_s(u_s)-\nu_s(u_s')+\big(\phi_s(u_s)-\phi_t(u_s')\big)x_s\Big)ds+\int_0^t\phi_s(u_s')y_sds\\
&+\int_0^t\Big(\chi_s(u_s)-\chi_s(u_s')+\big(\psi_s(u_s)-\psi_s(u_s')\big)x_s\Big)dB_s+\int_0^t\psi_s(u_s')y_sdB_s.
\end{align*}
Let
\begin{align*}
z_t=1-\int_0^t\phi_s(u_s')z_sds,
\end{align*}
and apply Ito's formula on $z_ty_t$ to get
\begin{align*}
z_ty_t&=\int_0^t\Big(\nu_s(u_s)-\nu_s(u_s')+\big(\phi_s(u_s)-\phi_t(u_s')\big)x_s\Big)z_sds\\
&+\int_0^t\Big(\chi_s(u_s)-\chi_s(u_s')+\big(\psi_s(u_s)-\psi_s(u_s')\big)x_s\Big)z_sdB_s\\
&+\int_0^tz_sy_s\psi_s(u'_s)dB_s\\
&=\int_0^t\Big(\nu_s(u_s)-\nu_s(u_s')+\big(\phi_s(u_s)-\phi_t(u_s')\big)x_s\Big)
\mathbb{I}_{\{u_s\neq u_s'\}}z_sds\\
&+\int_0^t\Big(\chi_s(u_s)-\chi_s(u_s')+\big(\psi_s(u_s)-\psi_s(u_s')\big)x_s\Big)
\mathbb{I}_{\{u_s\neq u_s'\}}z_sdB_s\\
&+\int_0^tz_sy_s\psi_s(u'_s)dB_s.
\end{align*}
 Then, using the Burkholder-Davis-Gundy and H\"older inequalities,
\begin{align*}
\E|zy|_T^{*,2p}&\leq K\bigg\{\bigg( \E\int_0^T\mathbb{I}_{\{u_t\neq u_t'\}}dt\bigg)^{1/2}+\int_0^T\E|z_sy_s|^{2p}ds\bigg\}\\
&=K\bigg\{\big(d(u,u')\big)^{1/2}+\int_0^T\E|z_sy_s|^{2p}ds\bigg\}.
\end{align*}
Thus, by Gronwalls inequality,
\begin{align*}
\E|zy|_T^{*,2p}\leq K\big(d(u,u')\big)^{1/2}.
\end{align*}
Finally, we conclude that
\begin{align*}
\E|x-x'|_T^{*,2p}=\E|z^{-1}zy|^{*,2p}\leq K\big(\E|zy|_T^{*,4p}\big)^{1/2}\leq K\big(d(u,u')\big)^{1/4},
\end{align*}
proving $(ii)$.

\medskip
\noindent Next, using the Lipschitz continuity of $h$ and $g$, we get
\begin{align*}
&|J(\delta_u)-J(\delta_{u'})|\\
\leq& \E\int_0^T\Big(|h(x_t,u_t)-h(x_t',u_t)|+|h(x_t',u_t)-h(x_t',u_t')|\Big)dt
+\E|g(x_T)-g(x_T')|\\
\leq& K\bigg\{\E\int_0^T\Big(|x_t-x_t'|+|h(x_t',u_t)-h(x_t',u_t')|\mathbb{I}_{\{u_s\neq u_s'\}}\Big)dt
+\E|x_T-x_T'|\bigg\}.
\end{align*}
 By $(ii)$ and the boundedness of $h$ we conclude that there exists $K>0$ such that
\begin{align*}
|J(\delta_u)-J(\delta_{u'})|&\leq K\big(d(u,u')\big)^{1/4}.
\end{align*}
\eop

\begin{lem} \label{adjointcont}
Let $u,u'\in\mathcal{U}$ along with the corresponding state processes $x_t,x_t'$ and adjoint processes $(p_t,q_t),
(P_t,Q_t),(p_t',q_t'),(P_t',Q_t')$ be given. Then it holds that
\begin{align}
\E\int_0^T\Big(|p_t-p_t'|^2+|q_t-q_t'|^2\Big)dt&\leq K\big(d(u,u')\big)^{1/8}, \label{adjointcont11}\\
\E\int_0^T\Big(|P_t-P_t'|^2+|Q_t-Q_t'|^2\Big)dt&\leq K\big(d(u,u')\big)^{1/8}, \label{adjointcont12}
\end{align}
for some constant $K>0$.
\end{lem}
\noindent\textbf{Proof.}
Let $(\bar{p}_t,\bar{q}_t)=
(p_t-p_t',q_t-q_t')$, which then satisfies the BSDE
\begin{align*}
\left\{ \begin{array}{ll}
d\bar{p}_t=&-\Big(\phi_t(u_t')\bar{p}_t+\psi_t(u_t')\bar{q}_t+\bar{\xi}_t\Big)dt+\bar{q}_tdB_t\medskip\\
\bar{p}_T=&g_x(x_T)-g_x(x_T'),
\end{array}\right.
\end{align*}
where
\begin{align*}
\bar{\xi}_t=\Big(\phi_t(u_t)-\phi_t(u_t')\Big)p_t+\Big(\psi_t(u_t)-\psi_t(u_t')\Big)q_t+h_x(x_t,u_t)-h_x(x_t',u_t').
\end{align*}
Further, let
\begin{align*}
y_t=\int_0^t\Big(\phi_s(u_s)y_s+|\bar{p}_s|\textrm{sgn}(\bar{p}_s)\Big)ds+\int_0^t\Big(\psi_s(u_s)y_s+|\bar{q}_s|
\textrm{sgn}(\bar{q}_s)\Big)dB_s,
\end{align*}
where $\textrm{sgn}(q_t)=\big(\textrm{sgn}(q_t^1),\ldots,\textrm{sgn}(q_t^d)\big)$ for the $d$-dimensional vector $q_t$. Now, by
defining
\begin{align*}
z_t=1-\int_0^t\phi_s(u_s)z_sds,
\end{align*}
and applying Ito's formula on $z_ty_t$, we get
\begin{align*}
z_ty_t=\int_0^tz_s|\bar{p}_s|\textrm{sgn}(\bar{p}_s)ds+\int_0^t\Big(z_sy_s\psi_s(u_s)+z_s|\bar{q}_s|
\textrm{sgn}(\bar{q}_s)\Big)dB_s.
\end{align*}
Note that $\E|z^{-1}|_T^{*,p}<\infty$ for any $p\geq 1$ by (A.3). Since $\psi_t$ is bounded, the Burkholder-Davis-Gundy and Gronwall inequalities yields
\begin{align*}
\E|zy|_T^{*,2p}\leq K\E\int_0^T|z_t|^{2p}\Big(|\bar{p}_t|^{2p}+|\bar{q}_t|^{2p}\Big)dt,
\end{align*}
for any $p\in\mathbb{N}$.
We can conclude using the Cauchy-Schwarz inequality that
\begin{align*}
\E|y|_T^{2p}=\E|z^{-1}zy|_T^{2p}\leq K\bigg(\E\int_0^T|z_t|^{4p}\Big(|\bar{p}_t|^{4p}+|\bar{q}_t|^{4p}\Big)dt\bigg)^{1/2}<\infty,
\end{align*}
for any $p\in\mathbb{N}$.
Moreover, by applying Ito's formula on $\bar{p}_ty_t$ and taking expectations, we obtain
\begin{align}
&\E\int_0^T\Big(\bar{p}_t|\bar{p}_t|\textrm{sgn}(\bar{p}_t)+\bar{q}_t|\bar{q}_t|\textrm{sgn}(\bar{q}_t)\Big)dt \notag\\
=& \E\bigg(\int_0^T\bar{\xi}_ty_tdt+\big(g_x(x_T)-g_x(x_T')\big)y_T\bigg) \label{adjointcont2}\\
\leq& K\bigg\{\bigg(\E\int_0^T|\bar{\xi}_t|^{2}dt\bigg)^{1/2}+\Big(\E|g_x(x_T)-g_x(x_T')|^2\Big)^{1/2}\bigg\}\Big(\E|y|_T^{*,2}\Big)^{1/2}.\notag
\end{align}
Using that the left hand side of (\ref{adjointcont2}) is equal to
$\E\int_0^T\big(|\bar{p}_t|^2+|\bar{q}_t|^2\big)dt$, we conclude
 that
\begin{align*}
\E\int_0^T\Big(|\bar{p}_t|^2+|\bar{q}_t|^2\Big)dt&\leq K\bigg\{\bigg(\E\int_0^T|\xi_t|^{2}dt\bigg)^{1/2}+\Big(\E|g_x(x_T)-g_x(x_T')|^2\Big)^{1/2}\bigg\}.
\end{align*}
As for the first term on the right hand side, we have by the Lipschitz continuity of $h_x$ that
\begin{align*}
&\E\int_0^T|h_x(x_t,u_t)-h_x(x_t',u_t')|^2dt\\
\leq &K\E\int_0^T\Big(|h_x(x_t,u_t)-h_x(x_t',u_t)|^2+
|h_x(x_t',u_t)-h_x(x_t',u_t')|^2\Big)dt\\
\leq&K\E\int_0^T\Big(|x_t-x_t'|^2+|h_x(x_t',u_t)-h_x(x_t',u_t')|^2\mathbb{I}_{\{u\neq u'\}}\Big)dt
\leq K\big(d(u,u')\big)^{1/4}.
\end{align*}
Moreover, by the estimate
\begin{align*}
&\E\int_0^T\Big(\Big|\big(\phi_t(u_t)-\phi_t(u_t')\big)p_t\Big|^2+\Big|\big(\psi_t(u_t)-\psi_t(u_t')\big)q_t\Big|^2\Big)dt\\
=&\E\int_0^T\Big(\Big|\big(\phi_t(u_t)-\phi_t(u_t')\big)p_t\Big|^2\mathbb{I}_{\{u\neq u'\}}+\Big|\big(\psi_t(u_t)-\psi_t(u_t')\big)q_t\Big|^2\mathbb{I}_{\{u\neq u'\}}\Big)dt\\
\leq&K\big(d(u,u')\big)^{1/4},
\end{align*}
we deduce
\begin{align*}
\bigg(\E\int_0^T|\bar{\xi}_t|^{2}dt\bigg)^{\frac{1}{2}}\leq K\big(d(u,u')\big)^{1/8}.
\end{align*}
Finally, by the Lipschitz continuity of $g_x$, we have the stimate
\begin{align*}
\Big(\E|g_x(x_T)-g_x(x_T')|^2\Big)^{\frac{1}{2}}\leq K\big(d(u,u')\big)^{1/8},
\end{align*}
which proves (\ref{adjointcont11}). Obviously, the same arguments prove (\ref{adjointcont12}).
\eop

\medskip
\noindent\textbf{Proof of Proposition \ref{nearoptimality}.}
Recall that by Lemma \ref{conttrajectory},~ $J:(\mathcal{U},d)\mapsto \mathbb{R}$ is continuous. Thus, by Ekeland's variational principle, see Ekeland \cite{ekeland} or Zhou \cite{zhou}, we can find a strict control $\tilde{u}_t$ such that
\begin{align*}
d(u,\tilde{u})\leq \epsilon^{2/3}
\end{align*}
and
\begin{align}
J(\delta_{\tilde{u}})\leq J(\delta_{u})+\epsilon^{1/3}d(u,\tilde{u}), \label{nearoptimal1}
\end{align}
for every $u\in\mathcal{U}$.
Next, we consider the perturbed control $\tilde{u}^\theta$:
\begin{align*}
\tilde{u}^\theta_t=\left\{ \begin{array}{ll}
v& \textrm{for}~ t\in [\tau,\tau+\theta]\\
\tilde{u}_t & \textrm{otherwise}.\end{array}\right.
\end{align*}
Thus, (\ref{nearoptimal1}) and the fact that $d(\tilde{u}^\theta,\tilde{u})\leq \theta$ implies
\begin{align*}
J(\delta_{\tilde{u}^\theta})-J(\delta_{\tilde{u}})\geq -\epsilon^{1/3}\theta.
\end{align*}
By Lemma \ref{varineq} the left hand side is equal to
\begin{align*}
&\E\int_0^T\Big(\mathcal{H}^{(\tilde{x}_t,\tilde{u}_t)}(t,\tilde{x}_t,\tilde{u}_t)-\mathcal{H}^{(\tilde{x}_t,\tilde{u}_t)}(t,\tilde{x}_t,\tilde{u}^\theta_t)\Big)dt+o(\theta).
\end{align*}
Since
\begin{align*}
\E\int_0^T\mathcal{H}^{(\tilde{x}_t,\tilde{u}_t)}(t,\tilde{x}_t,\tilde{u}^\theta_t)dt=&\E\int_0^T\Big(H\big(t,\tilde{x}_t,\tilde{u}^\theta_t,\tilde{p}_t,\tilde{q}_t-\tilde{P}_t(\chi_t(\tilde{u}_t)+\psi_t(\tilde{u}_t)\tilde{x}_t)\big)\\
 &-\frac{1}{2}\big(\chi_t(\tilde{u}^\theta_t)+\psi_t(\tilde{u}^\theta_t)\tilde{x}_t\big)\tilde{P}_t\big(\chi_t(\tilde{u}^\theta_t)+\psi_t(\tilde{u}^\theta_t)\tilde{x}_t\big)\Big)dt,
\end{align*}
and $\tilde{u}^\theta$ differs from $\tilde{u}$ only on $[\tau,\tau+\theta]$, we get
\begin{align*}
\E\int_\tau^{\tau+\theta}\Big(\mathcal{H}^{(\tilde{x}_t,\tilde{u}_t)}(t,\tilde{x}_t,\tilde{u}_t)-\mathcal{H}^{(\tilde{x}_t,\tilde{u}_t)}(t,\tilde{x}_t,v)\Big)dt+ o(\theta)\geq -\epsilon^{1/3}\theta.
\end{align*}
If we divide by $\theta$ and let $\theta\to0$, this yields
\begin{align}
\E\Big(\mathcal{H}^{(\tilde{x}_\tau,\tilde{u}_\tau)}(\tau,\tilde{x}_\tau,\tilde{u}_\tau)-\mathcal{H}^{(\tilde{x}_\tau,\tilde{u}_\tau)}(\tau,\tilde{x}_\tau,v)\Big)\geq -\epsilon^{1/3}.\label{nearoptimal5}
\end{align}
The next step is to replace $\tilde{u}$ with $u$, i.e we want to estimate
\begin{align}
&\E\int_0^T\Big(\mathcal{H}^{(\tilde{x}_t,\tilde{u}_t)}(t,\tilde{x}_t,\tilde{u})-\mathcal{H}^{(\tilde{x}_t,\tilde{u}_t)}(t,\tilde{x}_t,\tilde{u}^\theta_t)\Big)dt\notag\\
&-\E\int_0^T\Big(\mathcal{H}^{(x_t,u_t)}(t,x_t,u_t)-\mathcal{H}^{(x_t,u_t)}(t,x_t,u^\theta_t)\Big)dt. \label{nearoptimal2}
\end{align} We do this term by term. With $b$ and $\sigma$ as in (\ref{linearcoeff}), we may write (\ref{nearoptimal2}) as
\begin{align}
&\frac{1}{2}\E\int_0^T\Big(\sigma(t,\tilde{x}_t,\tilde{u}^\theta_t)-\sigma(t,\tilde{x}_t,\tilde{u}_t)\Big)\tilde{P}_t\Big(
\sigma(t,\tilde{x}_t,\tilde{u}^\theta_t)-\sigma(t,\tilde{x}_t,\tilde{u}_t)\Big)dt\notag\\
-&\frac{1}{2}\E\int_0^T\Big(\sigma(t,x_t,u^\theta_t)-\sigma(t,x_t,u_t)\Big)P_t\Big(\sigma(t,x_t,u^\theta_t)-\sigma(t,x_t,u_t)\Big)dt\notag\\
+&\E\int_0^T\tilde{q}_t\Big(b(t,\tilde{x}_t,\tilde{u}^\theta_t)-b(t,\tilde{x}_t,\tilde{u}_t)\Big)dt-\E\int_0^Tq_t\Big(b(t,x_t,u^\theta_t)-b(t,x_t,u_t)\Big)dt\notag\\
+&\E\int_0^T\tilde{P}_t\Big(\sigma(t,\tilde{x}_t,\tilde{u}^\theta_t)-\sigma(t,\tilde{x}_t,\tilde{u}_t)\Big)dt-\E\int_0^TP_t\Big(\sigma(t,x_t,u^\theta_t)-\sigma(t,x_t,u_t)\Big)dt\notag\\
+&\E\int_0^T\Big(h(t,\tilde{x}_t,\tilde{u}^\theta_t)-h(t,\tilde{x}_t,\tilde{u}_t)\Big)dt-\E\int_0^T\Big(h(t,x_t,u^\theta_t)-h(t,x_t,u_t)\Big)dt.\label{nearoptimal4}
\end{align}
We start by estimating the third line in the above expression.
\begin{align*}
&\E\int_0^T\tilde{q}_t\Big(b(t,\tilde{x}_t,\tilde{u}^\theta_t)-b(t,\tilde{x}_t,\tilde{u}_t)\Big)dt-\E\int_0^Tq_t\Big(b(t,x_t,u^\theta_t)-b(t,x_t,u_t)\Big)dt\\
=&\E\int_0^T\Big(\tilde{q}_t-q_t\Big)\Big(b(t,\tilde{x}_t,\tilde{u}^\theta_t)-b(t,\tilde{x}_t,\tilde{u}_t)\Big)dt\\
&+\E\int_0^Tq_t\Big(b(t,\tilde{x}_t,\tilde{u}^\theta_t)-b(t,x_t,\tilde{u}^\theta_t)\Big)dt
+\E\int_0^Tq_t\Big(b(t,x_t,u_t)-b(t,\tilde{x}_t,\tilde{u}_t)\Big)dt\\
=&I_1+I_2+I_3.
\end{align*}
We have the following estimate of the first term, using Lemma \ref{adjointcont} and the integrability of the components of $b(\cdot,\cdot,\cdot)$.
\begin{align*}
I_1\leq& \bigg(\E\int_0^T|\tilde{q}_t-q_t|^{2}dt\bigg)^{1/{2}}\bigg(\E\int_0^T|b(t,\tilde{x}_t,\tilde{u}^\theta_t)-b(t,\tilde{x}_t,\tilde{u}_t)|^2dt\bigg)^{1/2}\\
\leq&K\big(d(u,\tilde{u})\big)^{1/{16}}\leq K\epsilon^{1/{24}}.
\end{align*}
As for the second term, we get using Lemma \ref{conttrajectory}
\begin{align*}
I_2\leq& \bigg(\E\int_0^T|q_t|^2dt\bigg)^{1/2}\bigg(\E\int_0^T\big|\phi_t(u_t^\theta)\big[\tilde{x}_t-x_t\big]\big|^2\bigg)^{1/2}\\ \leq & K\bigg(\E\int_0^T\big|\tilde{x}_t-x_t\big|^4\bigg)^{1/4}\leq K\big(d(u,\tilde{u})\big)^{1/16}\leq K\epsilon^{1/24}.
\end{align*}
Further,
\begin{align*}
&I_3=\E\int_0^Tq_t\Big(b(t,\tilde{x}_t,u_t)-b(t,\tilde{x}_t,\tilde{u}_t)\Big)dt+\E\int_0^Tq_t\Big(b(t,x_t,u_t)-b(t,\tilde{x}_t,u_t)\Big)dt\\
&\leq \bigg(\E\int_0^T|q_t|^2\bigg)^{1/2}\bigg(\E\int_0^T\big|\nu_t(\tilde{u}_t)-\nu_t(u_t)+\Big(\phi_t(\tilde{u}_t)
-\phi_t(u_t)\Big)\tilde{x}_t\big|^2\mathbb{I}_{\{\tilde{u}\neq u\}}dt\bigg)^{1/2}\\
&+K\epsilon^{1/24}\\
&\leq K(\epsilon^{1/6}+\epsilon^{1/24}).
\end{align*}
We can conclude that
\begin{align*}
&\E\int_0^T\tilde{q}_t\Big(b(t,\tilde{x}_t,\tilde{u}^\theta_t)-b(t,\tilde{x}_t,\tilde{u}_t)\Big)dt-\E\int_0^Tq_t\Big(b(t,x_t,u^\theta_t)-b(t,x_t,u_t)\Big)dt\\
&\leq K\epsilon^{1/24}.
\end{align*}
By similar calculations we are able to get the same estimates for the other terms in (\ref{nearoptimal4}). Combining this with (\ref{nearoptimal5}), and since $v\in U$ is arbitrary, (\ref{Hinequality}) follows.
\eop

\medskip
\noindent\textbf{Proof of Lemma \ref{convergenceH}.}  
Let $(p^{(k)},q^{(k)}),(P^{(k)},Q^{(k)})$ and $(\hat{p},\hat{q}),(\hat{P},\hat{Q})$
 be the processes given by (\ref{adjoint11}) and (\ref{adjoint12}) corresponding to
$x_t^{(k)}$ and $\hat{x}_t$ respectively.  Define
\begin{align*}
\hat{z}_t&=1-\int_0^t\phi_s(\hat{\mu}_s)\hat{z}_sds,~\textrm{and}~z^{(k)}_t=1-\int_0^t\phi_s(u^{(k)}_s)z^{(k)}_sds.
\end{align*}
Further, we denote
$(\bar{p}_t,\bar{q}_t)=(\hat{z}_t\hat{p}_t-z^{(k)}_tp_t^{(k)},\hat{z}_t\hat{q}_t-z^{(k)}_tq_t^{(k)})$,
and simlarly for the second order adjoint processes (with $\phi$ replaced by
$2\phi$ above). The proof is carried out in three steps:
\begin{itemize}
\item[$(i)$] $\lim_{k\to\infty}\E\int_0^T\Big(\bar{p}_t^2+|\bar{q}_t|^2\Big)dt=0,$
\item[$(ii)$] $\lim_{k\to\infty}\E\int_0^T\Big(\bar{P}_t^2+|\bar{Q}_t|^2\Big)dt=0,$
\item[$(iii)$] $\lim_{k\to\infty}\E\int_0^T\mathcal{H}^{(x_t^{(k)},u_t^{(k)})}(t,x_t^{(k)},u_t^{(k)})dt=\E\int_0^T\mathcal{H}^{(\hat{x}_t,\hat{\mu}_t)}(t,\hat{x}_t,\hat{\mu}_t)dt.$
\end{itemize}

\medskip
We first prove ($i$). By applying It\^os formula on $\bar{p}_t^2$ we have
\begin{align*}
\bar{p}_t^2+\int_t^T|\bar{q}_s|^2ds
&=\big(\hat{z}_Tg_x(\hat{x}_T)-z^{(k)}_Tg_x(x_T^{(k)})\big)^2\\
&+2\int_t^T\bar{p}_s\Big(\hat{\Phi}_s-\Phi_s^{(k)}\Big)ds-2\int_t^T\bar{p}_s\bar{q}_sdB_s,
\end{align*}
where        
\begin{align*}
\hat{\Phi}_t=&\psi_t(\hat{\mu}_t)\hat{z}_t\hat{q}_t+\hat{z}_th_x(\hat{x}_t,\hat{\mu}_t)\\
\Phi^{(k)}_t=&\psi_t(u_t^{(k)})z^{(k)}_tq^{(k)}_t+z^{(k)}_th_x(x^{(k)}_t,u^{(k)}_t).
\end{align*}
Taking expectations and using Young's inequality on the second term on the right
hand side yields
\begin{align*}
\E\bar{p}_t^2+\E\int_t^T|\bar{q}_t|^2ds\leq&
\E\Big(\hat{z}_Tg_x(\hat{x}_T)-z_T^{(k)}g_x(x_T^{(k)})\Big)^2
+\alpha^2\E\int_t^T\bar{p}_s^2ds\\
&+\frac{1}{\alpha^2}\bigg(\E\int_t^T\Big(\hat{\Phi}_s-\Phi_s^{(k)}\Big)ds\bigg)^2.
\end{align*}
Expanding the last term:
\begin{align*}
&\E\int_t^T\Big(\hat{\Phi}_s-\Phi_s^{(k)}\Big)ds\\
=&\E\int_t^T\Big(\psi_s(\hat{\mu}_s)\hat{z}_s\hat{q}_s-\psi_s(u^{(k)}_s)z^{(k)}_sq_s^{(k)}+\hat{z}_sh_x(\hat{x}_s,\hat{\mu}_s)-z^{(k)}_sh_x(x_s^{(k)},u_s^{(k)})\Big)ds\\
\leq
&K\bigg\{\E\int_t^T|\psi_s(\hat{\mu}_s)-\psi_s(u^{(k)}_s)|^2|\hat{z}_s\hat{q}_s|^2ds+\E\int_t^T|\hat{z}_s\hat{q}_s-z^{(k)}_sq_s^{(k)}|^2|\psi_s(u^{(k)}_s)|^2ds\\
&+\E\int_t^T|\hat{z}_s-z^{(k)}_s||h_x(\hat{x}_s,\hat{\mu}_s)|ds+\E\int_t^T|z^{(k)}_s||h_x(\hat{x}_s,\hat{\mu}_s)-h_x(\hat{x}_s,u_s^{(k)})|\\
&+\E\int_t^T|z^{(k)}_s||h_x(\hat{x}_s,u_s^{(k)})-h_x(x_s^{(k)},u_s^{(k)})|ds\bigg\},
\end{align*}
where the first and fourth term converges to $0$ since $\delta_{u_t^{(k)}}dt\to
\hat{\mu}_t(du)dt$ $\P$-a.s. in $\mathcal{L}([0,T]\times U)$. Similarly, the third
term converges to $0$ since $\E|z^{(k)}-\hat{z}|_T^{*,2p}\to0$ for any $p\geq 1$. The
last term converges to $0$ since $h_x$ is
Lipschitz continuous and $\E|x^{(k)}-\hat{x}|_T^{*,2}\to0$. Inserting this expression into the above, and since $\psi$ is
bounded by some constant $C$, choose $\alpha$ such that $CK/\alpha^2<1$, it follows
by applying Gronwall's
inequality that
\begin{align*}
\E\bar{p}_t^2+\E\int_t^T|\bar{q}_t|^2ds\leq&
K\E\Big(\hat{z}_Tg_x(\hat{x}_T)-z_T^{(k)}g_x(x_T^{(k)})\Big)^2.
\end{align*}
Observing that $\E\big(g_x(\hat{x}_T)-g_x(x_T^{(k)})\big)^4\to 0$ by the Bounded Convergence
Theorem as well as \\$\E|z^{(k)}-\hat{z}|_T^{*,4}\to0$, the result follows.

\medskip
\noindent ($ii$) is proven with the same arguments.

\medskip
\noindent Now proving ($iii$) is straightforward by making use of ($i$), ($ii$) and the Chattering Lemma.
\eop

\end{document}